\documentclass{siamart190516}
\usepackage{amsmath}
\usepackage{graphicx}
\usepackage{float}
\usepackage{hyperref}
\usepackage{cleveref}
\usepackage{amssymb}
\usepackage[ruled,vlined,linesnumbered,algo2e]{algorithm2e}
\usepackage{dsfont}
\usepackage{pgfplots}
\usepackage{mathtools}
\usepackage{siunitx}
\usepackage[ruled,vlined]{algorithm2e}
\usepackage{filecontents}
\usepackage[none]{hyphenat}
\usepackage{verbatim}

\newcommand\smashtimes{\mathbin{\mathpalette\xhash\relax}}
\newcommand{\xhash}[2]{\ooalign{%
  $#1\xxhash{#1}{-45}$\cr
  $#1\xxhash{#1}{45}$\cr
  }%
}
\newcommand{\D}{\ensuremath{\mathrm{d}}}
\newcommand{\xxhash}[2]{\rotatebox[origin=c]{#2}{$#1\parallel$}}
\newcommand*\dif{\mathop{}\!\mathrm{d}}
\newcommand*\tset{\mathbf{t}}
\newcommand*\sset{\mathbf{s}}
\newcommand*\C{\mathbb{C}}

\newcommand*\Hmat{\mathcal{H}\text{-matrix}}
\DeclarePairedDelimiter\ceil{\lceil}{\rceil}

\theoremstyle{definition}
\newtheorem{example}[theorem]{Example}

\DeclareMathSymbol{\shortminus}{\mathbin}{AMSa}{"39}
\theoremstyle{remark}
\newtheorem{remark}[theorem]{Remark}
\DeclareMathOperator*{\argmax}{arg\,max}

\begin{document}
\sloppy
\title{Frequency extraction for BEM-matrices arising from the 3D scalar Helmholtz equation}
\headers{Frequency extraction 3D Helmholtz BEM-matrices}{}
\author{S. Dirckx, D. Huybrechs, K. Meerbergen}
\maketitle





\date{}


\newcommand*{\RED}[1]{{\color{red}#1}}
\newcommand*{\GREEN}[1]{{\color{green}#1}}

\begin{abstract}
\footnotesize
The discretisation of boundary integral equations for the scalar Helmholtz equation leads to large dense linear systems. Efficient boundary element methods (BEM), such as the fast multipole method (FMM) and $\Hmat$ based methods, focus on structured low-rank approximations of subblocks in these systems. It is known that the (numerical) ranks of these subblocks increase linearly with the wavenumber. We explore a data-sparse representation of BEM-matrices valid for a range of frequencies, based on extracting the known phase of the Green's function. Algebraically, this leads to a Hadamard product of a frequency matrix and an $\Hmat$. We show that the frequency dependence of this $\Hmat$ can be determined using a small number of frequency samples, even for geometrically complex three-dimensional scattering obstacles. We describe an efficient construction of the representation by combining adaptive cross approximation with adaptive rational approximation in the continuous frequency dimension. We show that our data-sparse representation allows to efficiently sample the full BEM-matrix at any given frequency, and as such it may be useful as part of an efficient sweeping routine.
\end{abstract}

\maketitle

\section{Introduction}
The scalar Helmholtz equation
\begin{equation}\label{eq:Helmholtz}
\nabla^2 \psi + \kappa^2 \psi = 0 \quad,\quad \text{on } \Omega\subset \mathbb{R}^3,\,\kappa\in\mathbb{R}
\end{equation}
models time-harmonic wave-propagation and
arises in the study of accoustic, electromagnetic and seismic waves. Both finite element methods (FEM) and boundary element methods (BEM) are used for the low- to mid-frequency regimes of $\kappa$. The FEM framework is attractive due to the sparseness of the resulting systems and the large body of existing methods, including eigenvalue computation, time marching and the fast computation of frequency response functions (FRF), all based on the sparseness of the systems. In addition the dependence of the system matrices on the wavenumber is simply quadratic. The boundary element method relies on the re-writing of the differential equation (\ref{eq:Helmholtz}) as a boundary integral equation (BIE) using the method of Green's functions \cite{McLean,SauterSchwab}. The discretised systems of the BIE formulation are dense and depend highly nonlinearly on the wavenumber $\kappa$. The advantage of BEM is  that only the boundary needs to be discretised, leading to a smaller system of equations having to be solved. Certainly for exterior problems (e.g. light scattering off a surface) and whenever only the boundary is of concern (e.g. seismic waves on the hull of a structure) this advantage is significant.

For moderate frequencies, the dense BEM matrices can be approximated by data-sparse matrices, either explicitly or implicitly. This results in, amongst others, the Hierarchical Matrix ($\Hmat$) framework \cite{Bebendorf,HMAT}, the Fast Multipole Method (FMM) \cite{Rohklin1,Rohklin3,Gumerov}, the panel clustering method \cite{Nowak} and wavelet based methods \cite{Utzinger,Huybrechs1}. These methods reduce the cost of assembling and solving the BEM-systems to near-linear. However, as the wavenumber $\kappa$ increases, so does the complexity of these methods, reducing their usefulness in the mid- to high-frequency range. In addition, the complicated wavenumber dependence of the BEM matrices is a largely open problem, the solution of which allows for new methods to be developed, including frequency sweeping and eigenvalue computations.

In this paper, we propose implicitly writing the Galerkin discretisation matrix of two main BEM operators as a Hadamard product of a frequency matrix with an $\Hmat$ i.e.
$$B(\kappa)=H(\kappa)\circ\hat{B}(\kappa)$$
with $H(\kappa)_{ij}:=\exp(\imath\kappa\, d_{ij})$, in which $d_{ij}$ is the distance between the centers of the $i$th test DOF and the $j$th ansatz DOF. We show that the matrix $\hat{B}$ has a significantly reduced assembly time compared with $B$ and allows for a convenient compact representation of the wavenumber dependence $\kappa\mapsto\hat{B}(\kappa)$, based on the (set-valued) AAA algorithm \cite{KarlSVAAA,TrefethenAAA,AAA}. This in turn sets up the use of Compact rational Krylov methods (CORK) for frequency sweeping \cite{CORK}.

This paper has two aims. The more important one is to introduce a compact representation for the wavenumber dependence of these BEM matrices that can be constructed in near-linear time as a function of the grid size. Secondly, we provide a method for the fast construction of BEM matrices arising from the BIE framework, whose data-sparseness is preserved at high frequencies. The Hadamard product is not easy to use in standard linear algebra algorithms. At this stage, we do not make any claims about the gain in complexity of matrix-vector product algorithms or iterative solvers, but only about the gain in the representation of the frequency dependent matrix.

The problem of increasing complexity for oscillatory kernels has been tackled in a number of ways, but the main focus has been that of \emph{directional} approximations. For instance, in \cite{Brandt}, a multilevel collocation scheme is outlined in which at each level the oscillatory part of the kernel is broken up into a number of directions, which increases with the coarsness of the level. This leads to a large constant appearing in the asymptotic complexity, which increases with the required accuracy. Other analytical schemes, such as the FMM methods in \cite{Messner} and \cite{BormDirectional} use Chebyshev interpolation on a modified kernel, obtained through multiplication by plane waves. These plane waves are defined at the level of interacting admissible clusters by the directions between cluster centers. As in \cite{Engquist}, they obtain specialised admissibility criteria that ensure low rank of admissible blocks. However, these admissibility criteria include an aperture requirement of $\mathcal{O}(1/\kappa)$. In the framework of $\mathcal{H}$-matrices, this means there are many small subblocks to be compressed and subsequently managed, resulting in high overhead and again a large fixed constant in the asymptotic complexity. In \cite{Brick}, a similar scheme with the same drawbacks is devised. An algebraic scheme based on this idea of plane wave multiplication is outlined in \cite{BebendorfDirectional}. Here, the same idea is proposed, only algebraically, using nested adaptive cross approximation rather than polynomial interpolation, constructing an $\mathcal{H}^2$-matrix, rather than an FMM scheme. In \cite{Butterfly}, a rank $r$ approximation of the phase function is exploited to construct a matrix-vector product through a multitude of $r$-dimensional non-uniform FFTs (NUFFT). This method is the closest to the one outlined in this paper, though we do not extract a low-rank approximant of the phase. In addition, the computational cost of NUFFTs quickly increases with the rank $r$. None of the methods outlined above tackle the problem of linearizing the wave number dependence of the BEM matrices.

\section{Boundary element methods and Hierarchical matrices}
Throughout this paper, $\mathbb{R}$ and $\mathbb{C}$ denote the real and complex number fields respectively. The symbol $\Omega$ denotes a bounded open domain with orientable (weakly) Lipschitz boundary $\partial\Omega$.
The main operators of interest in this paper are the \emph{single-layer potential} and \emph{double-layer potential} operator; these are operators on Sobolev spaces on the boundary $\partial\Omega$:
\begin{align*}
\mathcal{S}:H^{\frac{\shortminus1}{2}}(\partial\Omega)&\to H^{\frac{1}{2}}(\partial\Omega):\phi\mapsto\mathcal{S}\phi\\
\mathcal{D}:H^{\frac{1}{2}}(\partial\Omega)&\to H^{\frac{1}{2}}(\partial\Omega):\psi\mapsto\mathcal{D}\psi
\end{align*}
with (weakly) singular integral representations
\begin{align*}
\forall \phi\in L^{\infty}(\partial\Omega):(\mathcal{S}\phi)(\mathbf{y})&:=\int_{\partial\Omega}G(\mathbf{x},\mathbf{y})\phi(\mathbf{x})\dif S_{\mathbf{x}}\\
\forall \psi\in L^{\infty}(\partial\Omega):(\mathcal{D}\psi)(\mathbf{y})&:=\int_{\partial\Omega}\frac{\partial G}{\partial\mathbf{n}_{\mathbf{x}}}(\mathbf{x},\mathbf{y})\psi(\mathbf{x}) \dif S_{\mathbf{x}}
\end{align*}
in which
$$G(\mathbf{x},\mathbf{y}):=\frac{\exp(\imath\kappa \|\mathbf{x}\shortminus\mathbf{y}\|)}{4\pi\|\mathbf{x}\shortminus\mathbf{y}\|}$$
is the Green's function for the Helmholtz equation and $\mathbf{n}$ denotes the outward pointing unit normal vector field on $\partial\Omega$. The normal field $\mathbf{n}$ exists almost everywhere because $\partial\Omega$ was assumed to be (weakly) Lipschitz.
The operators $\mathcal{S}$ and $\mathcal{D}$ play an important role in the solution of the Helmholtz equation by means of integral equations, which we will illustrate with the example of the \emph{interior Dirichlet problem}:
\begin{definition}[Interior Dirichlet problem (IDP)]
Let $\Omega$ be an open bounded connected domain in $\mathbb{R}^3$ with boundary $\partial\Omega$. Let $\gamma_0 u$ denote the generalisation of the trace $u_{|\partial\Omega}$ to the Sobolev space $H^1(\Omega)$. The Dirichlet problem is to find $u\in H^1(\Omega)$ such that
\begin{alignat*}{2}
(\nabla^2+\kappa^2)u&=0 & &\text{ on } \Omega\\
\gamma_0u&=g_D & &\text{ on } \partial\Omega
\end{alignat*}
given Dirichlet data $g_D\in H^{\frac{1}{2}}(\partial\Omega)$. 
\end{definition}
The IDP can be reformulated as an integral equation as follows.
Let $\sigma(\mathbf{x})$ denote the surface angle of $\partial\Omega$ at $\mathbf{x}$, defined formally as 
\begin{equation}\label{eq:sigma}
\sigma(\mathbf{x})=\lim_{r\to 0}\frac{\mu(\partial B(\mathbf{x},r)\cap \Omega)}{\mu(\partial B(\mathbf{x},r))}
\end{equation}
with $\mu$ the standard surface measure on $\partial B(\mathbf{x},r)$.
Note that $\sigma=\frac{1}{2}$ almost everywhere, since we have assumed the unit normal field $\mathbf{n}$ to be well-defined almost everywhere on $\partial\Omega$.

Suppose $u$ is a solution to the IDP.
Let $u_N$ denote the associated Neumann data\footnote{ Note that since $u$ is an element of the Sobolev space $H^1(\Omega)$ this is not simply $\frac{\partial}{\partial\mathbf{n}} u$ but rather its generalisation to $H^{1}(\Omega)$ called the \emph{conormal trace} or \emph{conormal derivative} and is usually denoted as $u_{N}=\gamma_1 u$. See \cite{McLean} for a good introduction.}. Then (see \cite[\S3.4]{SauterSchwab})
\begin{equation*}
g_D=\shortminus \mathcal{S}u_{N} - ((\sigma-1)\mathds{1}+\mathcal{D})g_D
\end{equation*}
from which it follows that, for $\psi=u_N$
$$\mathcal{S}\psi=(\sigma\mathds{1}+\mathcal{D})g_D.$$
In addition, this integral formulation, whenever solvable, also defines a solution to the original IDP i.e. if $\psi$ solves the above integral equation, then it is the Neumann data of a solution of the IDP, which can then be reconstructed if needed (see \cite[\S7]{McLean}).
Similar integral formulations can be derived for the interior/exterior Dirichlet, interior/exterior Neumann and interior/exterior Robin problems. These are abbreviated to
\begin{equation}\label{eq:BEMGeneral2}
\mathcal{B}[\psi]=g\qquad\text{ on }\partial\Omega.
\end{equation}
For an overview of the possible systems, see \cite[\S3.4]{SauterSchwab}.

\subsection{Variational formulation}

The Boundary Element Method (BEM) seeks to solve the integral formulation (\ref{eq:BEMGeneral2}) using a variational approach. A simple introduction can be found in \cite{Pechstein}. A more thorough resource is \cite{SauterSchwab}. In such an approach, the solution $\psi$ is assumed to lie in an \emph{ansatzspace} $H_a=\text{span}\{\psi_j\}_j$ of distributions, and the right-hand side $g$ is assumed to lie in a \emph{testspace} $H_t=\text{span}\{\phi_i\}_i$ of distributions.\footnote{$\{\psi_j\}_j$ and $\{\phi_i\}_i$ are assumed to be bases.} Equation (\ref{eq:BEMGeneral2}) then discretises to 
$$\forall i:\langle\mathcal{B}[\textstyle\sum_ja_j\psi_j],\phi_i\rangle=\langle g,\phi_i\rangle=\int_{\partial\Omega}g(\mathbf{x})\bar{\phi}_i(\mathbf{x})\dif S_{\mathbf{x}}$$ 
in which we have adopted the convention that the inner product $\langle\cdot,\cdot\rangle$ is conjugate linear in the second argument. The discretisation of $\mathcal{B}$ and $g$ are correspondingly $B_{ij}=\langle\mathcal{B}\psi_j,\phi_i\rangle$ and $g_i=\langle g,\phi_i\rangle$. If we suppose $H_a$ and $H_t$ are finite-dimensional, of dimensions $n_a$ and $n_t$ respectively, then $B\in\mathbb{C}^{n_t\times n_a}$ and $\mathbf{g}\in\mathbb{C}^{n_t}$. Explicitly:
\begin{equation}\label{eq:discretizedSystem}
B\mathbf{a}=\begin{bmatrix}
\langle \mathcal{B}\psi_1,\phi_1 \rangle & \langle \mathcal{B}\psi_2,\phi_1 \rangle &\cdots &\langle \mathcal{B}\psi_{n_a},\phi_1 \rangle\\
\langle \mathcal{B}\psi_1,\phi_2 \rangle & \langle \mathcal{B}\psi_2,\phi_2 \rangle &\cdots &\langle \mathcal{B}\psi_{n_a},\phi_2 \rangle\\
\vdots & \vdots &\ddots &\vdots\\
\langle \mathcal{B}\psi_1,\phi_{n_t} \rangle & \langle \mathcal{B}\psi_2,\phi_{n_t} \rangle &\cdots &\langle \mathcal{B}\psi_{n_a},\phi_{n_t} \rangle\\
\end{bmatrix}
\begin{bmatrix}
a_1\\
a_2\\
\vdots\\
a_{n_a}
\end{bmatrix}=
\begin{bmatrix}
g_1\\
g_2\\
\vdots\\
g_{n_t}
\end{bmatrix}.
\end{equation}
When the integral operator is $\mathcal{S}$ or $\mathcal{D}$ instead of the abstract $\mathcal{B}$ we write $S$ and $D$ respectively for its discretisation.

In practice, the spaces $H_a$ and $H_t$ are built up from a mesh $\Gamma=\{\tau_n\}_{n=1}^{N}$ of the surface $\partial\Omega$. In this paper, $\Gamma$ is assumed to be triangular. Let $v(\Gamma)$ denote the set of vertices of $\Gamma$. Two common tactics exist:
\begin{itemize}
\item Galerkin discretisation: $\{\psi_j\}_j$ and $\{\phi_i\}_i$ are both sets of compactly supported functions. Typical examples are piecewise constant and piecewise linear continuous functions: 
\begin{align*}
S_0(\Gamma)&:=\{\chi_{\tau}\}_{\tau\in\Gamma}\\
S_1(\Gamma)&:=\{f_p\in C(\Gamma)|\forall \tau\in\Gamma: f_{p|\tau} \text{is linear and } \forall q\in v(\Gamma):f_p(q)=\delta_{pq}\}_{p\in v(\Gamma)}
\end{align*}
\item Collocation discretisation: $H_a$ is still some space of compactly supported functions (as in Galerkin), while $H_t=\{\delta_{p_j}\}_j$, where the collocation points $\{p_j\}_j$ are the centers of $\{\tau_j\}_{j}$.
\end{itemize}
Collocation, while faster to evaluate, has less desirable convergence and stability traits than Galerkin discretisation. Thus the focus in much of the literature and in this paper is on the Galerkin method. As in the finite element method (FEM), both the basis functions and their support are referred to by \emph{degrees-of-freedom} (DOFs). This is justified since the types of test function or ansatz function introduced above are uniquely defined by their support. We fix the following convention.
\begin{definition}
Throughout this paper, the index sets $I$ and $J$ are used to mean the set of test DOFs and ansatz DOFs respectively. Elements $i\in I$ and $j\in J$ are used to mean \textbf{at the same time}
\begin{itemize}
\item test DOFs and ansatz DOFs,
\item corresponding supports in $\Gamma$,
\item corresponding matrix entries in $B$.
\end{itemize}
This way, all of the following make sense for the set $\mathbf{t}\times\mathbf{s}\subseteq I\times J$:
\begin{itemize}
\item $\{\phi_i\}_{i\in \mathbf{t}}\subseteq H_t$, $\{\psi_j\}_{j\in \mathbf{s}}\subseteq H_a$,
\item $\forall i\in \mathbf{t}: i \subset \Gamma$,
$\forall j\in \mathbf{s}: j \subset \Gamma$
\item $B_{|\mathbf{t}\times\mathbf{s}}\in \mathbb{C}^{\tset\times\sset}$.
\end{itemize}
\end{definition}

\subsection{Hierarchical matrices}
The $\Hmat$ format and its arithmetic were introduced in \cite{HackbuschPart1}. They are an important tool in the solution of BEM systems. An $\Hmat$ is a dense matrix that is nevertheless data-sparse, owing to the fact that it contains many large low-rank subblocks, corresponding to pairs of so-called \emph{admissible} clusters. Standard works on the subject, on which this section is largely based, include \cite{HMAT} and \cite{Bebendorf}.\\
In contrast to the discretisations of finite element method (FEM) operators, the system matrices in the BEM formulation are not sparse, due to the fact that the single layer and double layer kernels are non-local. However, the Green's function $G(\mathbf{x},\mathbf{y})$ is \emph{separable} (also called \emph{degenerate}) for low frequencies in large parts of $I\times J$, meaning it can be well-approximated as 
$$G(\mathbf{x},\mathbf{y})\approx\sum_{i=1}^{r}G_1^{(i)}(\mathbf{x})G_2^{(i)}(\mathbf{y}).$$
Here $r$ is called the \emph{degree of degeneracy}. Clearly, if the above expansion holds for the block $\mathbf{t}\times\mathbf{s}$, then $G(\mathbf{t},\mathbf{s})$ and consequently $S_{|\mathbf{t}\times\mathbf{s}}$ are of rank (at most) $r$.

Typical separable expansions of $G$ of low degree (see \cite[\S3]{Bebendorf}, \cite[\S5]{Gumerov} and \cite{AminiProfit}) are only accurate when $\mathbf{x}$ and $\mathbf{y}$ are \emph{sufficiently separated}. This can be assured in a subblock $\mathbf{t}\times \mathbf{s} \subset I\times J$ by requiring
$$\eta \,\text{dist}(\xi_{\mathbf{t}},\mathbf{s})>\text{diam}(\mathbf{s})$$
with $\eta$ a parameter that encapsulates the requirement `sufficiently separated', and $\xi_{\mathbf{t}}$ the Chebyshev center of $\mathbf{t}$. The symmetric version of the above condition is
\begin{equation}\label{eq:admiss1}
\eta \,\text{dist}(\mathbf{t},\mathbf{s})>\max\{\text{diam}(\mathbf{t}),\text{diam}(\mathbf{s})\}.
\end{equation}
Often, $\eta = 2$ is used (for instance, it is the default $\eta$-value in \cite{HLibPro}). The condition (\ref{eq:admiss1}) is often relaxed to 
\begin{equation}\label{eq:admiss2}
\eta \,\text{dist}(\mathbf{t},\mathbf{s})>\min\{\text{diam}(\mathbf{t}),\text{diam}(\mathbf{s})\}.
\end{equation}
Typical kernels that are separable when condition $\ref{eq:admiss1}$ holds include the Laplace kernel $\frac{1}{\|\mathbf{x}-\mathbf{y}\|}$ and the euclidean distance kernel $\|\mathbf{x}-\mathbf{y}\|$.
These are not the only possible such conditions. Also, operations such as (exactly) calculating the diameter of clusters and the distance between clusters are computationally expensive. This motivates us to introduce a more general definition:
\begin{definition}[geometric admissibility condition]
A geometric admissibility condition for the DOF product set $I\times J$ and the kernel operator discretisation $B\in\mathbb{C}^{I\times J}$ is a boolean map
$$\text{Adm}:\mathcal{P}(I)\times \mathcal{P}(J)\to \{0,1\}$$
together with cluster distance and cluster diameter estimates 
\begin{align*}
\text{dist}^{\ast}&:\mathcal{P}(I)\times \mathcal{P}(J)\to \mathbb{R},\\
\text{diam}^{\ast}&:\mathcal{P}(I)\cup \mathcal{P}(J)\to \mathbb{R}.
\end{align*}
such that
\begin{itemize}
\item $\text{Adm}(\mathbf{t},\mathbf{s}) \Rightarrow B_{|\mathbf{t}\times\mathbf{s}}$ is (approximately) low-rank
\item $(\mathbf{t}_0\subseteq\mathbf{t}\land\mathbf{s}_0\subseteq\mathbf{s}\land \text{Adm}(\mathbf{t},\mathbf{s}))\Rightarrow\text{Adm}(\mathbf{t}_0,\mathbf{s}_0)$,
\item $(\text{Adm}(\mathbf{t},\mathbf{s})\land\text{dist}^{\ast}(\mathbf{t},\mathbf{s}')\geq \text{dist}^{\ast}(\mathbf{t},\mathbf{s})\land\text{diam}^{\ast}(\mathbf{s}')\leq \text{diam}^{\ast}(\mathbf{s}))\Rightarrow\text{Adm}(\mathbf{t},\mathbf{s}')$ and vice versa for $\mathbf{t}'$.
\end{itemize}
A block $\mathbf{t}\times\mathbf{s}$ such that $\text{Adm}(\mathbf{t},\mathbf{s})$ holds is called an admissible block. Given $\mathbf{t}\subset I$, any geometric admissibility condition divides $J$ into $\mathbf{t}$'s \emph{far-field} and its \emph{near-field}:
\begin{align*}
\text{Far}(\mathbf{t})&:=\{j\in J | \text{Adm}(\mathbf{t},\{j\})\},\\
\text{Near}(\mathbf{t})&:=\{j\in J | \lnot\text{Adm}(\mathbf{t},\{j\})\}.
\end{align*}
The near- and far-field of $\mathbf{s}\subset J$ are defined analogously. The far-field of $I\times J$ is the collection of all admissible blocks, i.e., $\{\tset\times\sset\in I \times J | \sset\in \text{Far}(\tset)\lor \tset\in \text{Far}(\sset)\}$.
\end{definition}
Typically, $\text{diam}^{\ast}$ and $\text{dist}^{\ast}$ are based on \emph{bounding boxes} for the clusters. These can be axis-aligned or cluster aligned. In the latter case, the bounding boxes are determined using principal component analysis (see \cite[\S1.4.1]{Bebendorf}).

The strategy to reduce $B\in\mathbb{C}^{I\times J}$ to a data-sparse format now boils down to determining an appropriate block partition $P_{I\times J}$ and subsequently approximating any admissible block in it by a low-rank matrix. Determining the partition $P_{I\times J}$ is done based on partitions $P_{I}$ and $P_{J}$ of $I$ and $J$, which in turn are based on the geometry of $\Gamma$ and a (bounding box based) binary subdivision strategy for the clusters in $\Gamma$:
$$s:\mathcal{P}(I)\to\mathcal{P}(I)\times \mathcal{P}(I):\mathbf{t}\mapsto(\mathbf{t}_0,\mathbf{t}_1)\, , \quad \text{ s.t. } \mathbf{t}=\mathbf{t}_0\cup\mathbf{t}_1 \text{ with } \mathbf{t}_0\cap\mathbf{t}_1=\emptyset.$$
An overview of division strategies can be found in \cite[\S1.4.1]{Bebendorf} and \cite[\S5.4]{HMAT}.

At level $0$, $P_I(0)$ is set to $I$. For level $l+1$ we say$P_I(l+1)$ is the collection of clusters resulting from the divisions of clusters in $P_I(l)$ according to the strategy $s$. This is repeated until a desired number of levels, say $L$, is reached. Formally, this is captured in the definition of the \emph{cluster tree} $\mathcal{T}_I$, a binary tree whose levels are precisely those described above. Note that its leaves are $P_I:=P_I(L)$. The cluster tree $\mathcal{T}_J$ is constructed in the same way.
From these two cluster trees, the \emph{block cluster tree} $\mathcal{T}_{I\times J}$ is constructed in the obvious way: its root is $I\times J$ and for each non-admissible node $\mathbf{t}^{(l)}_n\times \mathbf{s}^{(l)}_n$ at level $l$, its children are the elements of $s(\mathbf{t}^{(l)}_n)\times s(\mathbf{s}^{(l)}_n)$.\footnote{Sometimes, see \cite{HMAT}, blocks are divided regardless of admissibility. This is referred to as a \emph{level-conserving} block cluster tree. In that case, matrix block approximation is done not at the leaves of the block cluster tree, but rather at the highest admissible predecessor.} The leaves of the block cluster tree $\mathcal{T}_{I\times J}$ form the block cluster partition $P_{I\times J}$.
      
This finally enables us to formulate the definition of a \emph{Hierarchical matrix}.
\begin{definition}[Hierarchical matrix]\label{def:HMat}
A Hierarchical matrix ($\mathcal{H}$-matrix) is a matrix $M\in\mathbb{C}^{I\times J}$ together with a block cluster tree $\mathcal{T}_{I\times J}$ and a (geometric) admissibility criterion such that for all admissible blocks $\mathbf{t}\times\mathbf{s}\in P_{I\times J}$:
$$\exists X\in \mathbb{C}^{\mathbf{t}\times k}, Y\in \mathbb{C}^{\mathbf{s}\times k}: \qquad M_{|\mathbf{t}\times \mathbf{s}}=XY^T$$
with $k:=k(\mathbf{t},\mathbf{s})\ll\min\{\#\mathbf{t},\#\mathbf{s}\}$.
Given $B\in\mathbb{C}^{I\times J}$, we will denote its $\mathcal{H}$-matrix approximation by $B^{\mathcal{H}}$. When we want to stress the wavenumber $\kappa$, we will write $B(\kappa)$ and $B^{\mathcal{H}}(\kappa)$ instead.
\end{definition} 
In \cite[\S3]{Bebendorf}, it is shown that the discretisation of any  kernel operator with \emph{asymptotically smooth} kernel can be well-approximated by an $\mathcal{H}$-matrix, and that the kernels of $\mathcal{S}$, $\mathcal{D}$ are in this class for sufficiently low wave numbers.

\subsection{Adaptive Cross Approximation}
Low-rank approximation strategies based on explicit separable expansions are fast, but result in a sub-optimal rank. They require in-depth knowledge of the kernel, which is not always available. This is the reason why in many applications \emph{algebraic low-rank approximation schemes} are used. These take more time to construct, but do not require knowledge of the kernel function and result in a lower rank approximation. The most commonly used one, pioneered in \cite{Goreinov}, is known as the \emph{method of pseudoskeletons} or \emph{cross-approximation}. It can be made into an adaptive scheme, introduced in \cite{BebendorfACA}, called \emph{Adaptive Cross Approximation} (ACA).

\begin{algorithm2e}[ht]
\SetKwInOut{Input}{input}
\SetKwInOut{Output}{output}
\SetKw{Init}{init}{}{}
\SetAlgoLined
\Input{Matrix $M\in\mathbb{C}^{\tset\times \sset}$, tolerance tol}
\Output{Low-rank approx $M\approx XY^T$ }
\Init{$\epsilon := 1$, $k=0$}\\
 \While{$\epsilon >\text{tol}$}{
 select pivot $(i_k,j_k)\in \tset\times \sset$\\
 $\mathbf{x}_k:=M(:,j_k)$, $\mathbf{y}_k:=M(i_k,:)$\\
 $\mathbf{x}_k:=\mathbf{x}_k-\sum_{\mu<k}\mathbf{x}_{\mu}\mathbf{y}_{\mu}(j_k)$\\
 $\mathbf{y}_k:=\mathbf{y}_k-\sum_{\mu<k}\mathbf{x}_{\mu}(i_k)\mathbf{y}_{\mu}$\\
 $\mathbf{x}_k:=\frac{1}{\mathbf{x}_k(i_k)}\mathbf{x}_k$\\
  $X=[X,\mathbf{x}_k]$, $Y=[Y,\mathbf{y}_k]$\\
 $\epsilon := \|\mathbf{x}_k\|\|\mathbf{y}_k\|/(\|\mathbf{x}_0\|\|\mathbf{y}_0\|)$ \\
 $k:=k+1$
 }
 \caption{ACA}
 \label{alg:ACA}
\end{algorithm2e}
Intuitively, ACA sets an initial residue $R^{(0)}:=M$ and updates it to
\begin{equation}\label{eq:residue}
R^{(k+1)}:=R^{(k)}-\frac{1}{|R^{(k)}(i_k,j_k)|}R^{(k)}(:,j_k)R^{(k)}(i_k,:)
\end{equation}
until the update's norm is sufficiently small.
The matrix $R^{(k)}$ is never computed. Instead, the residue cross $(R^{(k)}(:,j_k),R^{(k)}(i_k,:))$ is calculated by generating the cross $(\mathbf{x}_k,\mathbf{y}_k^T)$ as in step \textbf{4}, and updating it in \textbf{5} and \textbf{6} which can be done in $\mathcal{O}(k(\#\tset+\#\sset))$ floating point operations (FLOP).

The most important step in ACA is the selection of the pivot. Ideally, $$(i_k,j_k)=\argmax_{(i,j)\in\tset\times\sset}|R^{(k)}(i,j)|,$$
but this would lead to a quadratic computational cost i.e. $\mathcal{O}(\#\tset\cdot\#\sset)$. One solution is ACA with \emph{partial pivoting} (see \cite[\S3.4]{Bebendorf}) or ACA+ (see e.g. \cite{HCA}). In ACA+, a reference cross $(\mathbf{x}_{\text{ref}},\mathbf{y}_{\text{ref}}^T)$ is first constructed, which is updated alongside the actual cross $(\mathbf{x}_k,\mathbf{y}_k^T)$. At every step, either $i_k$ or $j_k$ is set to  $\argmax_{i}|\mathbf{x}_{\text{ref}}(i)|$ or $\argmax_{j}|\mathbf{y}_{\text{ref}}(j)|$ respectively, depending on which of those attains the larger value. Then, the remaining pivot index ($j_k$ resp. $i_k$) is searched for over either $\mathbf{y}_k$ or $\mathbf{x}_k$.
\begin{remark}
The rank $k$ approximation $M^{(k)}=XY^{T}$ to $M$ may not be rank-optimal upon exit of the ACA+ algorithm. This can be remedied by \emph{algebraic recompression} (see \cite{Recompression}). Using the QR-decomposition of $X$ and $Y$ we write $XY^T$ as
\begin{align*}
XY^T &= Q_{X}R_XR_Y^TQ_Y^T\\
&=Q_{X}U\Sigma V^T Q_Y^T\\
&\approx Q_{X}U_{\epsilon}\Sigma_{\epsilon}V_{\epsilon}^TQ_Y^T\\
&=X_{\epsilon}\Sigma_{\epsilon}Y_{\epsilon}^T
\end{align*}
where $U\Sigma V^T$ and $U_{\epsilon}\Sigma_{\epsilon}V_{\epsilon}^T$ are the usual and the truncated singular value decomposition of the product $R_XR_Y^T$ respectively, with some prescribed truncation accuracy $\epsilon$.
The expression $X_{\epsilon}\Sigma_{\epsilon}Y_{\epsilon}^T$ constitutes a singular value decomposition, as $Q_X U_{\epsilon}$ and likewise $Q_Y V_{\epsilon}$ are orthogonal. This recompression can be done in $23k^3+6k^2(\# \mathbf{t}+\#\mathbf{s})$ floating point operations.
\end{remark}
Recall from equation \ref{eq:discretizedSystem} that for our purposes elements $M(i,j)$ of the admissible block $M$ are (theoretically) double integrals of the form 
\begin{equation}\label{eq:matrixElement}
M(i,j)=\langle\mathcal{B}\psi_i,\phi_j\rangle=\int_{i}\int_{j}b(\mathbf{x},\mathbf{y})\psi_i(\mathbf{y})\bar{\phi_j}(\mathbf{x})\dif S_{\mathbf{x}}\dif S_{\mathbf{y}}
\end{equation}
for a given kernel $b(\mathbf{x},\mathbf{y})$. In practice, the above is approximated using tensor quadrature via a Duffy transformation to a (2D) reference element (see \cite{GPGPUs}). Letting $C_{b}$ denote the number of floating point operations (FLOP) needed to evaluate the kernel $b$, and $q$ the order of Gaussian quadrature, it can then be shown that generating a single element of an admissible block requires $q^4(C_{b}+28)-1$ FLOP in the case of triangular elements and piecewise constant basis functions. Typical values for $q$ are $4$ or $5$. Due to the elementary functions involved, the FLOP count for kernel evaluation is dependent on the system and system settings. For $b=G$, on our machine and with our configuration, kernel evaluation was determined by high-precision timing to take $\sim 27$ FLOP. For simplicity, we say that $M(i,j)$ can be computed in $q^4\tilde{C}_b$ FLOP.
The complexity of ACA+ as described in \cite{HCA} is given by the following lemma:
\begin{lemma}\label{lem:costACA+}
Let the cost of generating an element $M(i,j)$ be given by $q^4\tilde{C}_b$ floating point operations (FLOP), counting multiplications, and additions. Suppose the ACA+ algorithm returns a rank $R$ matrix. Then the computational cost, say $c$, of the ACA+ algorithm, expressed in FLOP (counting additions, multiplications and other elementary operations\footnote{square root is counted as 5 FLOP}) is bounded by
$$ c < 2R(\#\tset+\#\sset)\left(q^4\tilde{C}_b+R+2\right)+R(2\max\{\#\tset,\#\sset\}+5).$$
Not included in this are single boolean comparisons and the integer updates of the form $k:=k+1$.
\end{lemma} 
\begin{proof}
This follows from careful counting, assuming that the worst case holds: the reference pivots are updated at each iteration.  
\end{proof}
A more useful and elegant complexity bound is to simply say that
$$c<2R(\#\tset+\#\sset)\left(q^4\tilde{C}_b+3R\right).$$
This holds because $\max\{\#\tset,\#\sset\}<(\#\tset+\#\sset)$ and because $2(\#\tset+\#\sset)>5$ for nontrivial blocks. In addition this assumes $2\leq R$, which is typically satisfied. Note that this complexity bound splits into two terms: one linear in both the block size and the rank, but which is of order 4 in $q$, and one term that is independent of quadrature and linear in the block size, but quadratic in the rank.
\section{Challenges of the $\mathcal{H}$-matrix format}
\subsection{Mid- and high-frequency regime}

Due to the oscillatory nature of the single- and double-layer potential kernels, admissible subblocks become dense in the limit $\kappa\to\infty$. The memory usage and consequently the time-complexity increase linearly with the wave number. This is true for both analytic and algebraic compression schemes.

 This is a huge challenge to using $\mathcal{H}$-matrix methods for the Helmholtz equation, in the mid- to high-frequency regime, and is one of the two main problems that we address in this paper. The main culprits of the increase in memory usage for the $\mathcal{H}$-matrix approximation to the operators $\mathcal{S}$ and $\mathcal{D}$ are the largest subblocks. For sufficiently small subblocks, a low-rank approximation is still possible. In fact, we can divide the set $P^{\text{Adm}}\subset P_{I\times J}$ of admissible subblocks into $P^{\text{large}}$ and $P^{\text{small}}$, based on the proportion
$\max\{\#\mathbf{t},\#\mathbf{s}\}/\kappa$. The idea (see \cite[\S10.4]{HMAT} and \cite{H2Hackbusch}) is then to write
$$B = B^{\mathcal{H}} + B^{\mathcal{H}^2},$$
where $B^{\mathcal{H}}$ is an $\mathcal{H}$-matrix with $P^{\text{small}}$ as block cluster partition, that is zero on the blocks in $P^{\text{large}}$. The matrix $B^{\mathcal{H}^2}$ is zero everywhere \emph{except} in the blocks in $P^{\text{large}}$, where a low-rank approximation is constructed using a separable expansion (e.g. Amini-Profit multipole expansion (see \cite{AminiProfit})). This defines an $\mathcal{H}^2$-structure, as introduced in \cite{H2SourceBorm} and \cite{H2Source}.

While this has a better (asymptotic) performance, it has serious drawbacks. Firstly, even for 2D problems, the memory usage of $B$ still increases linearly with the wave number $\kappa$. Secondly, the format of $B^{\mathcal{H}^2}$ does not lend itself well to efficient matrix operations. Thirdly, the evaluation of the multipole expansion coefficients can become numerically unstable. Fourthly, the implementation of the $\mathcal{H}^2$ format and methods is much more involved than that of the standard $\mathcal{H}$-matrix format.

The above approach is similar to high-frequency FMM already proposed in 1989 in \cite{Rohklin3}. Here, the high-frequency kernel is treated with a separable analytical expansion as well. The application of the discrete SLP and DLP to a vector can be done quickly and stably using `diagonal translation operators', as in \cite{Rohklin2}. Modern implementations (e.g. \cite{HFFMM} and \cite{HFFMM2}) have a complexity of $\mathcal{O}(N\log N)$.

Both approaches outlined above have good asymptotic behavior, but with large overhead and a large constant.
One goal of our paper therefore is simple: To overcome the highly oscillatory behaviour of the kernel at high wavenumbers without resorting to more involved data formats or possibly unstable expansions while maintaining a near-linear order.

\subsection{Compact representation of the wavenumber dependence}
The most important goal of this paper is to obtain a compact representation of the wavenumber dependence $\kappa\mapsto B(\kappa)$. This is a highly nontrivial objective, as illustrated by example \ref{ex:dependenceExtracted}. In this example, the SLP matrix $S$ and its $\mathcal{H}$-matrix approximant were investigated at a single index of an admissible block. The low-rank compression was achieved using ACA+ with an error tolerance of $1e\text{-}5$. The example was run using the open-source package $\mathcal{H}-\text{lib}^{\text{pro}}$ (see \cite{HLibPro}) using fifth order quadrature. All other options were set to the default. The figure shows two things. The $\Hmat$ format (at least in its most basic form) is ill-equipped to deal with high frequency behaviour, as evidenced by the increasing error in this regime. At the same time $\kappa\mapsto B^{\mathcal{H}}(\kappa)$ is not a good approximation for $\kappa\mapsto B(\kappa)$. Furthermore, the figure shows that $\kappa\mapsto S_{ij}(\kappa)$ approximately contains only one frequency, whereas $\kappa\mapsto S_{ij}^{\mathcal{H}}(\kappa)$ clearly does not. In other words, the  numerical approximation using ACA+ introduces significant numerical noise, which makes the approximation \emph{qualitatively} undesirable. The wave number dependence approximation can perhaps be improved by careful parameter tuning or by using specialised low-rank approximation schemes. However, by the Nyquist-Shannon sampling theorem, the top two figures illustrate that this task can only be achieved using a large number of samples. This motivates us to explore an alternative approach in which this problem is entirely avoided.
\section{Frequency extraction of the kernel function}

The central observation of this paper is the following. Suppose $i\in I$ and $j\in J$ are in each other's farfield. Suppose, in addition, that $\xi_i$ and $\xi_j$ are their (Chebyshev) centers. Then, with $S$ the discretisation of the single layer potential operator $\mathcal{S}$, it holds that
\begin{align}\label{eq:splitting}
S_{ij}&=\int_{i}\int_{j} G(\mathbf{x},\mathbf{y})\bar{\phi}_i(\mathbf{x})\psi_j(\mathbf{y}) \dif S_{\mathbf{x}} \dif S_{\mathbf{y}}\nonumber\\
&=\frac{1}{4\pi}\int_{i}\int_{j} \frac{\exp(\imath \kappa\|\mathbf{x}-\mathbf{y}\|)}{\|\mathbf{x}-\mathbf{y}\|}\bar{\phi}_i(\mathbf{x})\psi_j(\mathbf{y})\dif S_{\mathbf{x}} \dif S_{\mathbf{y}}\nonumber\\
&=\frac{\exp(\imath \kappa\|\xi_i-\xi_j\|)}{4\pi}\int_{i}\int_{j} \frac{\exp(\imath \kappa(\|\mathbf{x}-\mathbf{y}\|-\|\xi_i-\xi_j\|))}{\|\mathbf{x}-\mathbf{y}\|}\bar{\phi}_i(\mathbf{x})\psi_j(\mathbf{y})\dif S_{\mathbf{x}} \dif S_{\mathbf{y}}.
\end{align} 
This constitutes a splitting of $S_{ij}$ into an oscillatory part -- the prefactor in \eqref{eq:splitting} -- and a separable part. The following theorem is proven in Appendix \ref{appendix:A}.
\begin{theorem}\label{thm:LRindep}
Let $\tset\subset I$ and $\sset\subset J$ be clusters such that $\text{dist}(\tset,\sset)>h\sqrt{5}/2$ where $h$ is the diameter of the largest DOF in $\tset\cup\sset$ and both $\kappa h <1/2$ and $\kappa\text{dist}(\tset,\sset)>1$ are assumed. Suppose in addition that $\tset\times\sset$ is admissible in the sense of \ref{eq:admiss1} with $\eta=1/2$. Let $H(\kappa)\in\mathbb{C}^{\tset\times\sset}$ be given by $H(\kappa)_{ij}:=\exp(\imath \kappa \|\xi_i-\eta_j\|)$, where $\xi_i$ and $\eta_j$ are the centers of the DOFs $i$ and $j$. Then for any fixed $0<\epsilon<1$, the numerical $\epsilon$-rank of $\bar{H}(\kappa)\circ S_{|\tset\times\sset}$ satisfies
$$\text{rank}_{\epsilon}(\bar{H}(\kappa)\circ S_{|\tset\sset})\leq\mathcal{O}\left(\frac{\log^6\epsilon}{\log^2\kappa h}\right)$$
\end{theorem}
The interpretation of this is that the rank of $\bar{H}(\kappa)\circ S_{\sset\times\tset}$ depends only \emph{weakly} on the wave number. If the clusters are kept the same (as sets in $\mathbb{R}^3$) and $\kappa$ is increased while keeping $\kappa h$ constant (i.e. the meshing of the clusters \emph{is} refined), the complexity does not increase linearly with $\kappa$. The requirement $\kappa h<1/2$, while sometimes violated in practice, is typical of theoretical derivations involving discretised BEM matrices. The above theorem is no longer effective when $\kappa$ is increased independently of $h$. Note that we have only assumed $\tset\times\sset$ to be admissible in a weaker sense of \ref{eq:admiss1}. We did not need to design specialised admissibility criteria.
This leads us to the following definition:
\begin{definition}[Frequency extraction]
Let $B(\kappa)$ be the Galerkin discretisation of either the single- or double-layer potential for the Helmholtz equation at wavenumber $\kappa$. Then writing $B(\kappa)$ as the Hadamard product
$$B(\kappa)=H(\kappa)\circ \hat{B}(\kappa)$$
with
\begin{equation}
 H_{ij}(\kappa):=
\begin{cases}
\exp(\imath\kappa\|\xi_i-\xi_j\|),&i\in \text{far}(j)\\
1,&\text{else}
\end{cases}
\end{equation}
is called frequency extraction. Matrix $\hat{B}$ is defined as the extracted discretisation matrix. 
\end{definition}
The advantage of the above representation should be clear. When using ACA, we avoid many of the costly quadrature based evaluations of the original matrix, since $\hat{B}$ has much lower rank in its admissible blocks. Frequency extracion uses only the DOF centers and elementary functions, so it greatly reduces the overall cost, when applied in one of two ways. Firstly it is possible to approximately generate BEM-matrices $B(\kappa)$ in a data sparse format for higher frequencies, using the fact that $B=H\circ \hat{B}$, where $\hat{B}$ is a hierarchical matrix. This method of generating BEM matrices will be referred to as `one-shot frequency extraction'. Secondly, as we will show in section \ref{sec:Compact}, extracting the oscillatory behaviour of the kernel permits the construction of a compact representation of the wavenumber dependence $\kappa\mapsto \hat{B}(\kappa)$.
\subsection{One-shot frequency extraction}
In the one-shot method, rather than constructing an $\Hmat$ approximation of the discretisation $B$, we construct one for the extracted $\hat{B}$. This means that the matrix $B$ cannot be computed explicitly, as the Hadamard product
$$B=H\circ \hat{B}$$
would result in a quadratic cost. The main advantage of this method is the reduction in storage cost. Certainly in the case of very large or coupled systems, this is advantageous. One of the important drawbacks is the introduction of the Hadamard product, which is not easy to exploit in standard linear algebra algorithms. However, the special structure of $H(\kappa)$ suggest that a fast matrix-vector product based on a fast matrix-vector product with $\hat{B}$ exist. In particular, NUFFT-based (see \cite{NUFFT}) or Butterfly algorithm-based (see \cite{Butterfly}) approaches seem promising. The use of the Hadamard product in BEM solvers remains an open problem.

\subsection{Frequency extraction for the compact representation of the wavenumber dependence}\label{sec:Compact}
Consider the following example:
\begin{example}\label{ex:dependenceExtracted}
Let $\Gamma$ be the standard triangular mesh for the 3D sphere $S^2$ with $N=131072$ faces. Let $S,S^{\mathcal{H}},D,D^{\mathcal{H}}\in\mathbb{C}^{I\times J}$ denote the discretisations of the single- and double-layer potential and the $\Hmat$ approximation thereof, using $I=J=S_0(\Gamma)$ and ACA+. Then, with $(i,j)=(48640,6418)$, we obtain Figure ~\ref{fig:dep3}  for the wavenumber dependencies. Similar observations can be made for the double layer potential.
\begin{figure}
\center
\begin{tikzpicture}
\begin{axis}[width=.4\linewidth,
  xlabel=$\kappa\cdot\text{diam}(\Gamma)$,
  ylabel=$\text{im}({S_{ij},\hat{S}_{ij}})$,
  ylabel style={font=\small},
  grid,
  grid style={dotted}]
\addplot+ [mark size=1pt] table {A_ij_im_bf.dat};
\addplot+ [mark size=1pt] table {B_ij_im_bf.dat};
\end{axis}
\end{tikzpicture}
\begin{tikzpicture}
\begin{axis}[width=.4\linewidth,
  xlabel=$\kappa\cdot\text{diam}(\Gamma)$,
  ylabel=$\text{re}({S_{ij},\hat{S}_{ij}})$,
  ylabel style={font=\small},
  grid,
  grid style={dotted}]
\addplot+ [mark size=1pt] table {A_ij_real_bf.dat};
\addplot+ [mark size=1pt] table {B_ij_real_bf.dat};
\end{axis}
\end{tikzpicture}\\
\begin{tikzpicture}
\begin{axis}[width=.4\linewidth,
  xlabel=$\kappa\cdot\text{diam}(\Gamma)$,
  ylabel=$\text{im}({S^{\mathcal{H}}_{ij},\hat{S}^{\mathcal{H}}_{ij}})$,
  ylabel style={font=\small},
  grid,
  grid style={dotted}]
\addplot+ [mark size=1pt] table {A_ij_im.dat};
\addplot+ [mark size=1pt] table {B_ij_im.dat};
\end{axis}
\end{tikzpicture}
\begin{tikzpicture}
\begin{axis}[width=.4\linewidth,
  xlabel=$\kappa\cdot\text{diam}(\Gamma)$,
  ylabel=$\text{re}({S^{\mathcal{H}}_{ij},\hat{S}^{\mathcal{H}}_{ij}})$,
  ylabel style={font=\small},
  grid,
  grid style={dotted}]
\addplot+ [mark size=1pt] table {A_ij_real.dat};
\addplot+ [mark size=1pt] table {B_ij_real.dat};
\end{axis}
\end{tikzpicture}

\caption{Top: original BEM operator discretisation $S_{ij}(\kappa)$ (blue) together with its extracted discretisation $\hat{S}_{ij}(\kappa)$ (red). Bottom: $\Hmat$ approximations $\hat{S}_{ij}^{\mathcal{H}}(\kappa)$ and $\hat{S}^{\mathcal{H}}$. All as a function of the dimensionless wave number $\kappa\cdot\text{diam}(\Gamma)$.}
\label{fig:dep3}
\end{figure}
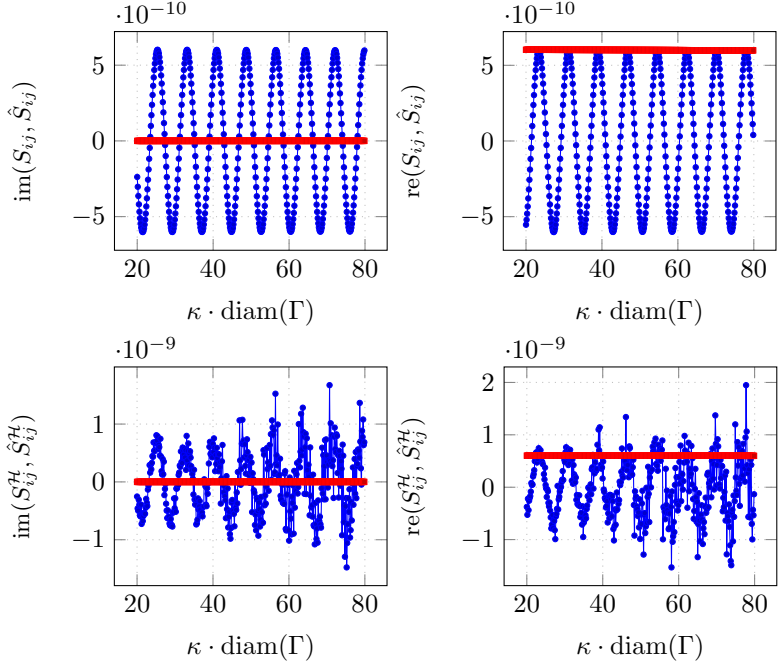
\end{example}
We observe that the extracted kernel operators are far less oscillatory and far better approximated by an $\Hmat$ than their non-extracted counterpart. In particular, for fixed $(i,j)$, this shows that $\hat{S}_{ij}(\kappa)$, $\hat{S}_{ij}^{\mathcal{H}}(\kappa)$, $\hat{D}_{ij}(\kappa)$ and $\hat{D}_{ij}^{\mathcal{H}}(\kappa)$ can be well approximated by polynomials or rational functions. Our preference is for rational functions, since rational approximation is a generic tool for which many algorithms exist. In addition, rational functions have proven to be easily exploitable in eigenvalue computation algorithms and fast frequency sweeps (see \cite{CORK}).
The construction of a compact wave number dependence representation, which requires additionally that the wave number dependence is low-rank over $\tset\times\sset$, is the subject of section \ref{sec:6}.

\section{AAA-ACA and rational approximation}\label{sec:6}
We now model the wavenumber dependence $\kappa\mapsto\hat{B}(\kappa)$ block-wise by interpreting $\kappa\mapsto\hat{B}_{|\tset\times\sset}(\kappa)$ as an element of the tensor product space $\mathbb{C}^{\tset\times\sset }\otimes C^1([a,b])$, for all admissible $\tset\times\sset$. Here $[a,b]$ is the frequency range of $\kappa$, with $0<a<b<\infty$.

\begin{definition}[Tensor approximation]
Let $\mathcal{E} \subset C^{1}([a,b])$ (see remark \ref{rem:proximinal}). Consider $\hat{\mathcal{B}}_{\tset\times\sset}\in\mathbb{C}^{\tset\times\sset }\otimes C^1([a,b])$ defined as
$$\hat{\mathcal{B}}_{\tset\times\sset}:[a,b] \to \mathbb{C}^{\tset\times\sset}:\kappa \mapsto \hat{B}_{|\tset\times\sset}(\kappa)$$
for admissible $\tset\times\sset$.
The approximation 
$$\hat{\mathcal{B}}_{\tset\times\sset}\approx \hat{\mathcal{B}}^{\otimes}_{\tset\times\sset}=\sum_{k=1}^{R_T}M_k \otimes f_k, \quad \text{ with }\quad\forall k: M_k\in\mathbb{C}^{\tset\times\sset},\,f_k\in \mathcal{E}$$
with low $R_T$ is called the low tensor rank approximation of $\hat{\mathcal{B}}_{|\tset\times\sset}$\footnote{$(M_k \otimes f_k)(i,j,\kappa)=M(i,j)\cdot f(\kappa)$}. The rank $R_T$ is called the tensor rank. The matrix function $\hat{\mathcal{B}}^{\otimes}$ defined by
$$\hat{\mathcal{B}}^{\otimes}_{|\tset\times\sset}:= \hat{\mathcal{B}}^{\otimes}_{\tset\times\sset}.$$
is called the tensor approximant of the matrix function $\kappa\mapsto\hat{B}(\kappa)$.
\end{definition}
\begin{remark}\label{rem:proximinal}
The set $\mathcal{E}\subset C([a,b])$ should be suited for efficient and effective function approximation. This can be any proximinal set, including, but not limited to, polynomials of fixed degree, rational functions of fixed degree with poles outside $[a,b]$, exponential sums and splines. Our choice is to set $\mathcal{E}$ to the set of rational functions with poles outside $[a,b]$.
\end{remark}
An approximate projection onto the set of rational functions of fixed degree with poles outside of $[a,b]$ is given by the AAA algorithm from \cite{AAA}. This rational approximation algorithm, characterised by a set of nodes $S$, maps a function $f$ to $\textbf{aaa}(f,S)\in\mathcal{E}$. Using this, a tensor ACA scheme can be built. This is given in algorithm \ref{alg:AAA-ACA}.
\begin{algorithm2e}[ht]
\label{alg:AAA-ACA}
\caption{AAA-ACA}
\KwData{ AAA algorithm $f\mapsto \texttt{aaa}(f,S) $, interpolation points $S$, ACA routine \textbf{aca}}
\SetKwInOut{Input}{input}
\SetKwInOut{Output}{output}
\SetKw{Init}{init}{}{}
\SetAlgoLined
\Input{Matrix function $\hat{\mathcal{B}}_{|\tset\times\sset}\in\mathbb{C}^{\tset\times \sset}\otimes C^1([a,b])$, tolerance \textbf{tol}, maximal rank $R_{\max}$}
\Output{Low tensor rank approx $\hat{\mathcal{B}}_{|\tset\times\sset}\approx \sum_{k=1}^{R_T} M_k\otimes f_k$ }
\Init{$\epsilon := \infty$, $r=1$}\\
 \While{$\epsilon >\textbf{tol}$ \emph{\&} $k<R_{\max}$}{
 select pivot $((i_k,j_k),\kappa_k)\in (\tset\times\sset)\times S$\hspace*{\fill}($\ast$)\label{line:pivot}\\
 $f_k:\texttt{aaa}([\kappa\mapsto \hat{\mathcal{B}}_{i_kj_k}(\kappa)])$\\
 $M_k:=\mathbf{aca}(\hat{\mathcal{B}}_{|\tset\times\sset}(\kappa_k))$\label{line:acaMk}\\
 $f_k:=f_k-\sum_{\mu<k}(M_\mu)_{i_k,j_k}f_{\mu}$\\
 $M_k:=M_k-\sum_{\mu<r}f_{\mu}(\kappa_k)M_{\mu}$ \hspace*{\fill}($\ast$)\label{line:Mk}\\
 $f_k:=f_k/f_k(\kappa_k)$\\
 $\epsilon := \|f_k\|\|M_k\|_{F}/(\|f_0\|\|M_0\|_F)$ \hspace*{\fill}($\ast$)\label{line:Err}\\
 $k:=k+1$
 }
\end{algorithm2e}
In fact, algorithm \ref{alg:AAA-ACA} has precisely the same form as algorithm \ref{alg:ACA}, only tensorial, with the exception of line \ref{line:acaMk}. Setting $M_{k}:=\hat{\mathcal{B}}_{|\tset\times\sset}(\kappa_k)$ would immediately lead to a quadratic cost $\mathcal{O}(\#\tset\cdot\#\sset)$. This is remedied by setting instead $M_{k}:=\mathbf{aca}(\hat{\mathcal{B}}_{|\tset\times\sset}(\kappa_k))$, using any linear time ACA routine $\mathbf{aca}$ (e.g. ACA+). Understand that this means that $M_k = X_kY_k^T$ and could in principle be stored as such.
However, implementing the tensor approximation routine like in algorithm~\ref{alg:AAA-ACA} would still lead to a quadratic cost. The main reasons for this are as follows.
\begin{itemize}
\item In step \textbf{\ref{line:Mk}} $\hat{\mathcal{B}}_{|\tset\times\sset}(\kappa_k)$ might be well-approximated by a low-rank matrix, but
$$M_k-\sum_{\mu<k}f_{\mu}(\kappa_k)M_{\mu}$$
is in general not low-rank. We need to introduce an alternative data-sparse and computationally beneficial representation of $M_k$.
\item The selection of $(i_k,j_k)$ in line \textbf{\ref{line:pivot}} is over $\tset\times\sset$, which leads to a quadratic cost, no matter the strategy. Therefore the search space for $(i_k,j_k)$ needs to be restricted.
\item Calculating the Frobenius norm $\|M_k\|_F$ in step \textbf{\ref{line:Err}} is expensive. An error estimate is therefore needed.
\end{itemize}
\subsection{Linear time AAA-ACA} 
This section outlines the implementation of a linear time version of Algorithm~\ref{alg:AAA-ACA}. Firstly, we will outline how the matrix $M_r$ can be constructed implicitly and evaluated.
\begin{proposition}\label{prop:coeffs}
In the notation from Algorithm~\ref{alg:AAA-ACA}, suppose for $k\in\{1,\ldots,R_T\}$
$$\textbf{aca}(\hat{\mathcal{B}}(\kappa_{k}))=X_{k}Y_{k}^T.$$
Then
$$M_{k}=\sum_{\nu=1}^{k}c_{\nu,k}X_{\nu}Y_{\nu}^T$$
with $c_{k,k}=1$ and
$$c_{\nu,k}:=-\sum_{\mu=\nu}^{k-1}f_{\mu}(\kappa_{k})c_{\nu,\mu}$$
for all $\nu <k$.
\end{proposition} 
\begin{proof}
For $k=1$ the proposition is obvious, as
$$M_1=\textbf{aca}(\hat{\mathcal{B}}(\kappa_{1}))=X_1Y_1^T,$$
Suppose the proposition is true for $k-1$. For $k$ we then have
\begin{align*}
M_k&=X_kY_k^T-\sum_{\mu=1}^{k-1}f_{\mu}(\kappa_k)M_{\mu}\\
&=X_kY_k^T-\sum_{\mu=1}^{k-1}f_{\mu}(\kappa_k)\left(\sum_{\nu=1}^{\mu}c_{\nu,\mu}X_{\nu}Y_{\nu}^T\right)\\
&=X_kY_k^T-\sum_{\nu=1}^{k-1}\left(\sum_{\mu=\nu}^{k-1}f_{\mu}(\kappa_k)c_{\nu,\mu}\right)X_{\nu}Y_{\nu}^T.
\end{align*}
And so the proposition follows by induction.
\end{proof}
We can interpret the coefficients up to $k$ as a coefficient matrix $C_k$ defined by
\begin{equation}
\label{eqn:coeffs}
C_{k} := \begin{bmatrix}
1 & c_{1,2} & c_{1,3}&\hdots & c_{1,k}\\
& 1 & c_{2,3} & \hdots &c_{2,k}\\
&& \ddots& \ddots & \vdots\\
&&& 1 & c_{k\shortminus 1,k}\\
& & & &1 
\end{bmatrix}.
\end{equation}
Then, with $\mathbf{f}(\kappa_{k+1}):=[f_1(\kappa_{k+1}),\ldots,f_k(\kappa_{k+1})]^T$ we have
\begin{equation}
\label{eqn:coeffsUpdate}
C_{k+1} := \begin{bmatrix}
C_{k}&\shortminus C_k\mathbf{f}(\kappa_{k+1})\\
\mathbf{0}^T & 1
\end{bmatrix}.
\end{equation}
Note that the coefficients $\{c_{\nu,k}\}_{\nu\leq k}$ can be calculated using only $\frac{(k-1)k}{2}$ multiplications and $\frac{(k-2)(k-1)}{2}$ additions, for a total of $(k-1)^2$ floating point operations (FLOP). Note also that this proposition makes it unnecessary to store the dense matrices $\{M_k\}_{k=1}^{R_T}$ calculated in Algorithm~\ref{alg:AAA-ACA}. Rather, we store the set $\{C_{R_T},\{X_k,Y_k\}_{k=1}^{R_T},\{f_k\}_{k=1}^{R_T}\}$. The following definition formalises this. Efficient reconstruction from such a representation will be treated in section~\ref{sec:reconstruction}.
\begin{definition}\label{def:compact}
Let $\hat{\mathcal{B}}_{|\tset\times\sset}\approx \sum_{k=1}^{R_T} M_k\otimes f_{k}$ with
$$M_{k}=\sum_{\nu=1}^{k}c_{\nu,k}X_{\nu}Y_{\nu}^T$$
as in proposition \ref{prop:coeffs}. Let $C_{R_T}\in\C^{R_T\times R_T}$ be as in equation \ref{eqn:coeffs}. Then we define the corresponding \emph{compact representation} of $\hat{\mathcal{B}}_{\otimes|\tset\times\sset}$ as $\mathcal{C}_{\tset\times\sset} := \{C_{R_T},\{X_k,Y_k\}_{k=1}^{R_T},\{f_k\}_{k=1}^{R_T}\}$.\\
The collection $\{\mathcal{C}_{\tset\times\sset}\}_{\tset\times\sset\in\text{Adm}(\Gamma)}$ will be called the $\emph{compact representation}$ of $\hat{\mathcal{B}}^{\otimes}$.
\end{definition} 

Secondly, we construct an approach to the pivot search that relies only on a subset of $\tset\times\sset$. It suffices here to observe that adaptive cross-approximation can be tweaked such that it also outputs the chosen pivots. This leads to the following definition:
\begin{definition}[\textbf{aca}-skeleton]
Let \textbf{aca} be an adaptive cross-approximation scheme applied to $M_k\in\mathbb{C}^{\tset\times\sset}$. Let $S_k^{\tset}\subset\tset$ and $S_k^{\sset}\subset\sset$ denote the pivots selected by $\textbf{aca}$ in the construction of the low-rank approximant to $M_{\tset\times\sset}$. Then the set $$S_k^{\tset}\smashtimes S_k^{\sset}:=(S_k^{\tset}\times\sset)\cup(\tset\times S_k^{\sset})$$ is called the \textbf{aca}-skeleton of $M$. We write
$$(X_k,Y_k,S^{\tset}_k\smashtimes S^{\sset}_k)=\textbf{aca}(M_k)$$
to emphasise the pivots also being output.
\end{definition}
The approach is simple: restrict the search for $(i_k,j_k)$ in step \textbf{\ref{line:pivot}} to the \textbf{aca}-skeleton $(S_k^{\tset}\smashtimes S_k^{\sset})$. Clearly, this results in a cost of $\mathcal{O}(R_k(\#\tset+\#\sset))$, with $R_k$ the rank of the ACA-approximant of $M(\kappa_k)$.\\
Lastly, the Frobenius norm $\|M_k\|_F$ in step \textbf{\ref{line:Err}} is replaced by a sampled estimate:
\begin{definition}\label{def:errEst}
Let $M_k$ be as in Algorithm~\ref{alg:AAA-ACA}. We define the norm estimator $\|M_k\|_{F,m}$ as 
$$\|M_k\|_{F,m}^2:=\frac{(\#\tset\cdot\#\sset)}{m}\|M_{k|E_{m}}\|_F^2$$ with $E_m\subset\tset\times\sset$ a random subset of cardinality $m$, drawn from a uniform distribution without repetition.
\end{definition}
Obviously, if $m=\#\tset\cdot\#\sset$, the above estimate is exact. The goal is to find a low $m$ such that$|\,\|M_k\|_{F,m}-\|M_k\|_{F}|$
is sufficiently small. The choice we settled on, which gives good results, is to take
\begin{equation*}
m=R_k(\#\tset+\#\sset).
\end{equation*}
The AAA-ACA algorithm can now be reformulated as Algorithm~\ref{alg:linear AAA-ACA}.\\
\begin{algorithm2e}[ht]
\label{alg:linear AAA-ACA}
\caption{Linear time AAA-ACA}
\KwData{ AAA algorithm $f\mapsto \texttt{aaa}(f,S) $, interpolation points $S$, ACA routine \textbf{aca}}
\SetKwInOut{Input}{input}
\SetKwInOut{Output}{output}
\SetKw{Init}{init}{}{}
\SetAlgoLined
\Input{Matrix function $\hat{\mathcal{B}}\in\mathbb{C}^{\tset\times \sset}\otimes C^1(\tset\times\sset)$, tolerance tol, maximal tensor rank $R_{\max}$}
\Output{Compact representation $\{C_{R_T},\{X_k,Y_k\}_{k=1}^{R_T}\},\{f_k\}_{k=1}^{R_T}\}$ }
\Init{$\epsilon := \infty$, initialise $S^{\tset}_0\subset\tset$ and $S^{\sset}_0\subset\sset$, $k=1$}\\
 \While{$\epsilon >\text{tol}$ \emph{\&} $k<R_{\max}$}{
 Set $((i_k,j_k),\kappa_k):=$\newline
 $\hspace*{-5pt}\argmax_{\tiny{((i,j),\kappa)\in(S_{k-1}^{\tset}\smashtimes S_{k-1}^{\sset})\times S}}\hat{\mathcal{B}}_{ij}(\kappa)-\sum_{\mu<k-1}(\sum_{\nu\leq\mu}c_{\mu,\nu}(X_{\nu}Y_{\nu}^T)_{ij})f_{\mu}(\kappa)$\label{line:pivot2}\\
 $f_r:\texttt{aaa}([\kappa\mapsto \hat{\mathcal{B}}_{i_kj_k}(\kappa)])$\label{line:aaa}\\
 $(X_k,Y_k,S_k^{\tset}\smashtimes S_k^{\sset}):=\textbf{aca}(\hat{\mathcal{B}}_{\tset\times\sset}(\kappa_k))$\label{line:acaInFinal}\\
 calculate $\{c_{k,\nu}\}_{\nu=1}^k$ as in proposition \ref{prop:coeffs}\\
 $f_k:=f_k-\sum_{\mu<k}(M_{\mu})_{i_kj_k}f_{\mu}=f_k-\sum_{\mu<k}(\sum_{\nu\leq\mu}c_{\mu,\nu}(X_{\nu}Y_{\nu}^T)_{i_kj_k})f_{\mu}$\label{line:fk}\\
  $f_k:=f_k/f_k(\kappa_k)$\label{line:normalize}\\
 $\epsilon := \|f_k\|\|M_k\|_{F,m}/(\|f_1\|\|M_1\|_F)$\label{line:errInAAA}\\
 $k:=k+1$ 
 }
\end{algorithm2e}
Here the functions $f_k$ are stored and passed as vectors in $\mathbb{C}^{\#S}$, as the values on $S$ fully define $f_k$. Using lemma \ref{lem:costACA+}, careful counting and the $\textbf{aaa}$ complexity estimate shown in \cite{AAA}, we obtain the following:
\begin{proposition}
Let $R_T$ be the tensor rank attained by Algorithm~\ref{alg:linear AAA-ACA}. Suppose in addition that $R$ is the maximal rank attained by the $\textbf{aca}$ step and $n$ is the maximal number of nodes needed for interpolation in the \textbf{aaa} step. Then, if ACA+ is used as the $\textbf{aca}$ compression scheme and $\|f_k\|$ is computed using quadrature in $c_f\cdot (\#S)$ FLOP, the cost $c$ in FLOP of Algorithm~\ref{alg:linear AAA-ACA} is bounded by
$$c<\frac{2}{3}R_T^3K(R+1)+R_T^2K(1+2\alpha+2R)+R_TK(\beta+6R)$$
with $K=R(\#\tset+\#\sset)$, $\alpha=\frac{\#S}{K}$ and $\beta = \alpha(n^3c_{\textbf{aaa}}+c_f+1)+(3+\alpha)q^4\tilde{C}_{b}+3$, where $c_{\textbf{aaa}}$ is a constant such that the cost of the $\textbf{aaa}$ algorithm is bounded by $c_{\textbf{aaa}}\cdot(\# S)\cdot n^3$.
\end{proposition}
\begin{proof}
The proof of this is simple. At step $k$, the cost of index selection is bounded by $R(\#\tset+\#\sset)(q^4\tilde{C}_b+2R(1+\cdots+k)+k)$ FLOP. The cost of the $\textbf{aaa}$ and $\textbf{aca}$ step are bound respectively by $c_{\textbf{aaa}}\cdot(\#S)\cdot n^3$ and $2R(\#\tset+\#\sset)(q^4\tilde{C}_b+3R)$ FLOP, per \cite{AAA} and lemma~\ref{lem:costACA+}. The coefficient computation requires $(k-1)^2$ FLOP, by equation~\ref{eqn:coeffsUpdate}.
Line~\ref{line:fk}, \ref{line:normalize} and \ref{line:errInAAA} require respectively $(\#S)(q^4\tilde{C}_b+2k)+2R(1+\cdots+k)$, $(\#S)$ and $c_{f}(\#S)+2R(\#\tset+\#\sset)(1+kR)+7$ FLOP. Here we used that the square root requires about 5 FLOP. Then, using
$$\sum_{k=1}^{R_T}
(1+\cdots+k) = \frac{R_T(R_T+1)(R_T+2)}{6}\leq\frac{1}{2}R_T^3$$
and
$$\sum_{k=1}^{R_T}(k-1)^2 = \frac{(R_T-1)R_T(2R_T-1)}{6}<\frac{1}{3}R_T^3$$
for $2\leq R_T$, the proof readily follows under the assumption that $\#\tset,\#\sset \geq 2$, i.e. the block considered is nontrivial.
\end{proof}
Observe that this shows that Algorithm~\ref{alg:linear AAA-ACA} is linear in the block size.
\begin{remark}
Algorithm \ref{alg:linear AAA-ACA} is similar to the algorithms introduced in \cite{TuckerACA}, but differs in some key ways. Firstly, the intermediate matrices $M_k$ (see Algorithm~\ref{alg:AAA-ACA}) are not assumed to be numerically low-rank. In \cite{TuckerACA}, Algorithm 2 step (3), a low rank approximation of $M_k$ is made at each iteration. Computing these $M_k$ is relatively costly. Our Algorithm~\ref{alg:linear AAA-ACA}, through the coefficient matrix $C$, avoids the computation of the intermediate $M_k$ entirely. Only for the error computation are some elements of these $M_k$ used. Secondly, the error heuristic and pivot selection heuristic that we use are simpler, while still performing well. Thirdly, our method is more general in that it does not assume that the $\kappa$-direction is given in discretised form.
\end{remark}
\subsection{Reconstruction from the compact representation}\label{sec:reconstruction}
In this subsection we outline an efficient way to, from a given  tensor approximation $\hat{\mathcal{B}}^{\otimes}$, reconstruct an $\Hmat$ $\hat{\mathcal{B}}^{\mathcal{H}}(\kappa)$ for any given $\kappa$. In particular we provide a way to, from the representation $\{C_{R_T},\{X_k,Y_k\}_{k=1}^{R_T},\{f_k\}_{k=1}^{R_T}\}$ for a block $\hat{\mathcal{B}}^{\otimes}_{|\tset\times\sset}$, construct a low-rank approximation to $\hat{B}(\kappa)_{|\tset\times\sset}$. We provide estimates for the complexity of our method. We use the following:
\begin{lemma}\label{lem:reconstruction}
Suppose $\kappa$ is a given wave number. Let $$\hat{\mathcal{B}}^{\otimes}_{|\tset\times\sset}=\sum_{k=1}^{R_T}M_{k}\otimes f_{k}$$
such that, as in proposition \ref{prop:coeffs},
$$M_{k} = \sum_{\nu=1}^{k}c_{\nu,k}X_{\nu}Y_{\nu}^T.$$ Then using the notation from equation \ref{eqn:coeffs}, and $\mathbf{f}=[f_1(\kappa),\ldots,f_{R_T}(\kappa)]^T$ let
$$\mathbf{m}:=C_{R_T}\mathbf{f}.$$
It now holds that
$$\hat{B}^{\otimes}_{|\tset\times\sset}(\kappa) = \sum_{k=1}^{R_T} \mathbf{m}(k)X_{k}Y_{k}^T$$
\end{lemma}
\begin{proof}
This is a direct consequence of proposition \ref{prop:coeffs}.
\end{proof} 
Lemma \ref{lem:reconstruction} provides a convenient way to evaluate rows and columns of $\hat{B}^{\otimes}_{|\tset\times\sset}(\kappa)$ efficiently. This means that adaptive cross-approximation schemes can be used to construct a low-rank approximation to a block $\hat{B}^{\otimes}_{|\tset\times\sset}(\kappa)$. We have implemented a specialised ACA+ algorithm for this, which is just a simple re-implementation of ACA+ where $\hat{B}^{\otimes}_{|\tset\times\sset}(i,j)$ is evaluated not through quadrature, but through the procedure outlined in lemma~\ref{lem:reconstruction}. Lemma~\ref{lem:costACACrossList} is then simple to see, following lemma~\ref{lem:costACA+} and using the fact that $C_{R_T}\mathbf{f}$ can be computed in less than $R_T^2$ FLOP.
\begin{lemma}\label{lem:costACACrossList}
Suppose the rank of $\hat{\mathcal{B}}^{\otimes}_{|\tset\times\sset}(\kappa)$ found by ACA+ is $\overline{R}$ and that the maximal rank of $\{(X_k,Y_k)\}_{k=1}^{R_T}$ is $R$. Then the cost $c$ of the specialised ACA+ in FLOP is bounded by
$$c<2\overline{R}(\#\tset+\#\sset)\left(2RR_T+3\overline{R}\right)+2R_T^2.$$
\end{lemma}

\section{Numerical results: One-shot}
In this section we present our numerical results for the one-shot problem. All these experiments were run on a Linux x86\_64 machine with 8 Intel Core i7-6820HQ 2.70 GHz, 8192 KB cache processors each with 4 cores. We used the open-source package $\mathcal{H}-\text{lib}^{\text{pro}}$ (see \cite{HLibPro}) in which clustering, admissibility criteria, quadrature rules and low-rank approximation come already implemented.

\subsection{Shapes used in the experiments}
We chose four shapes to be used in our experiment: the sphere $S^1$, the NETGEN crankshaft (see \cite{netgen}), the Stanford bunny (see \cite{Stanford}) and the NASA Almond (see \cite{almond}). These four shapes were meshed, resulting in the (triangular) meshes seen below:
\begin{figure}[ht]
\center
\includegraphics[height=3cm,width=.3\linewidth]{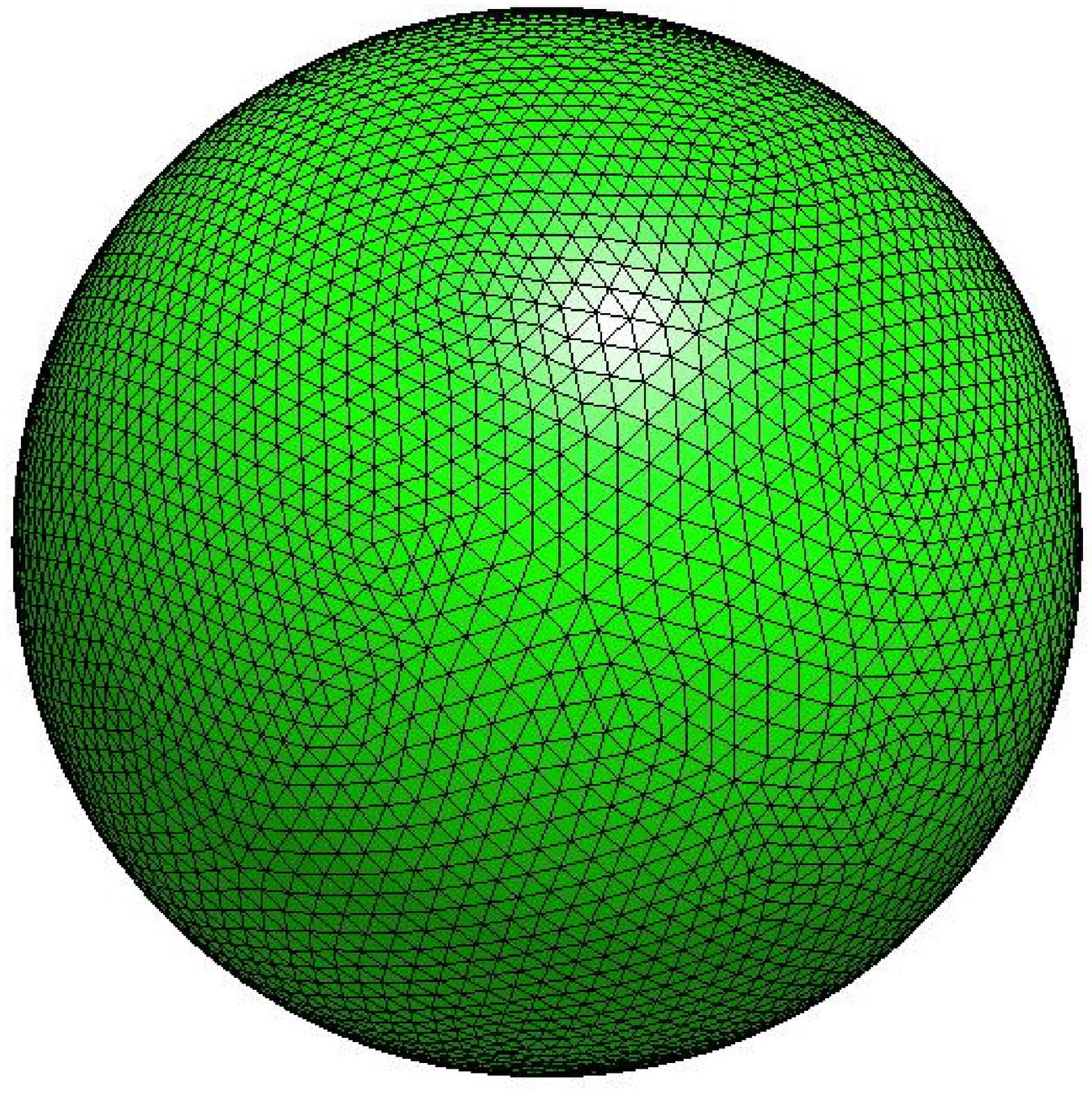}
\includegraphics[height=3cm,width=.3\linewidth]{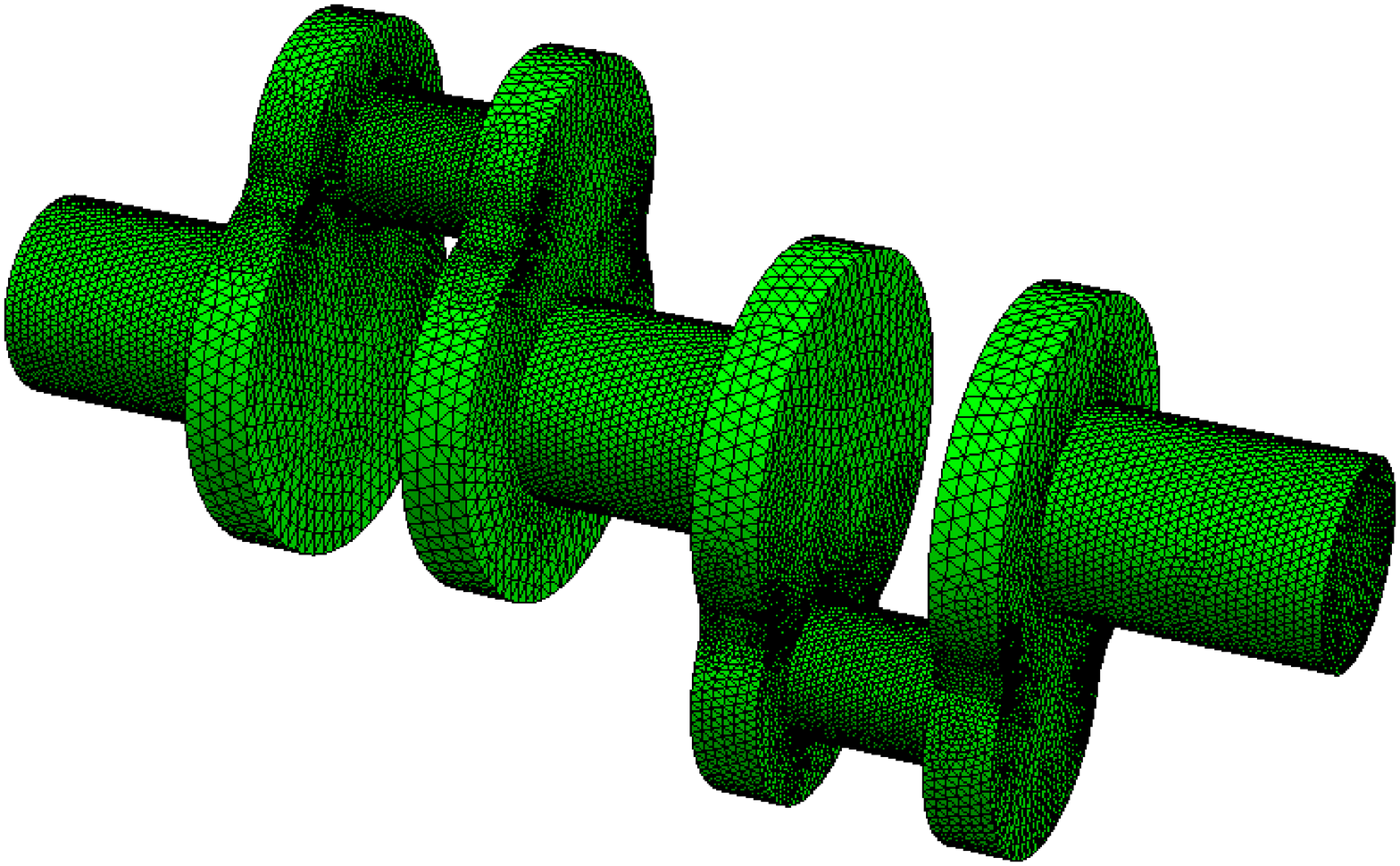}

\includegraphics[height=3cm,width=.3\linewidth]{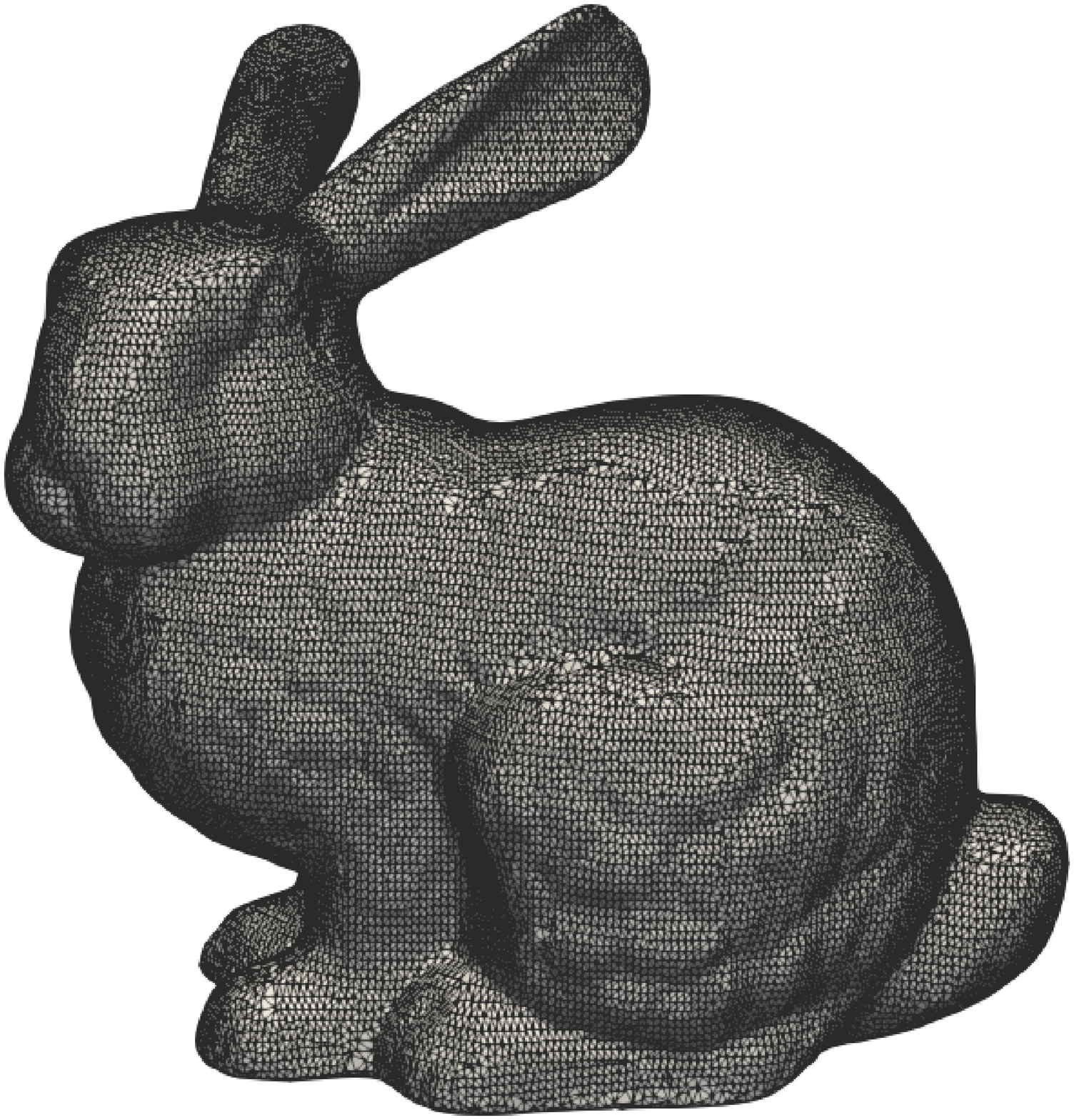}
\includegraphics[height=3cm,width=.3\linewidth]{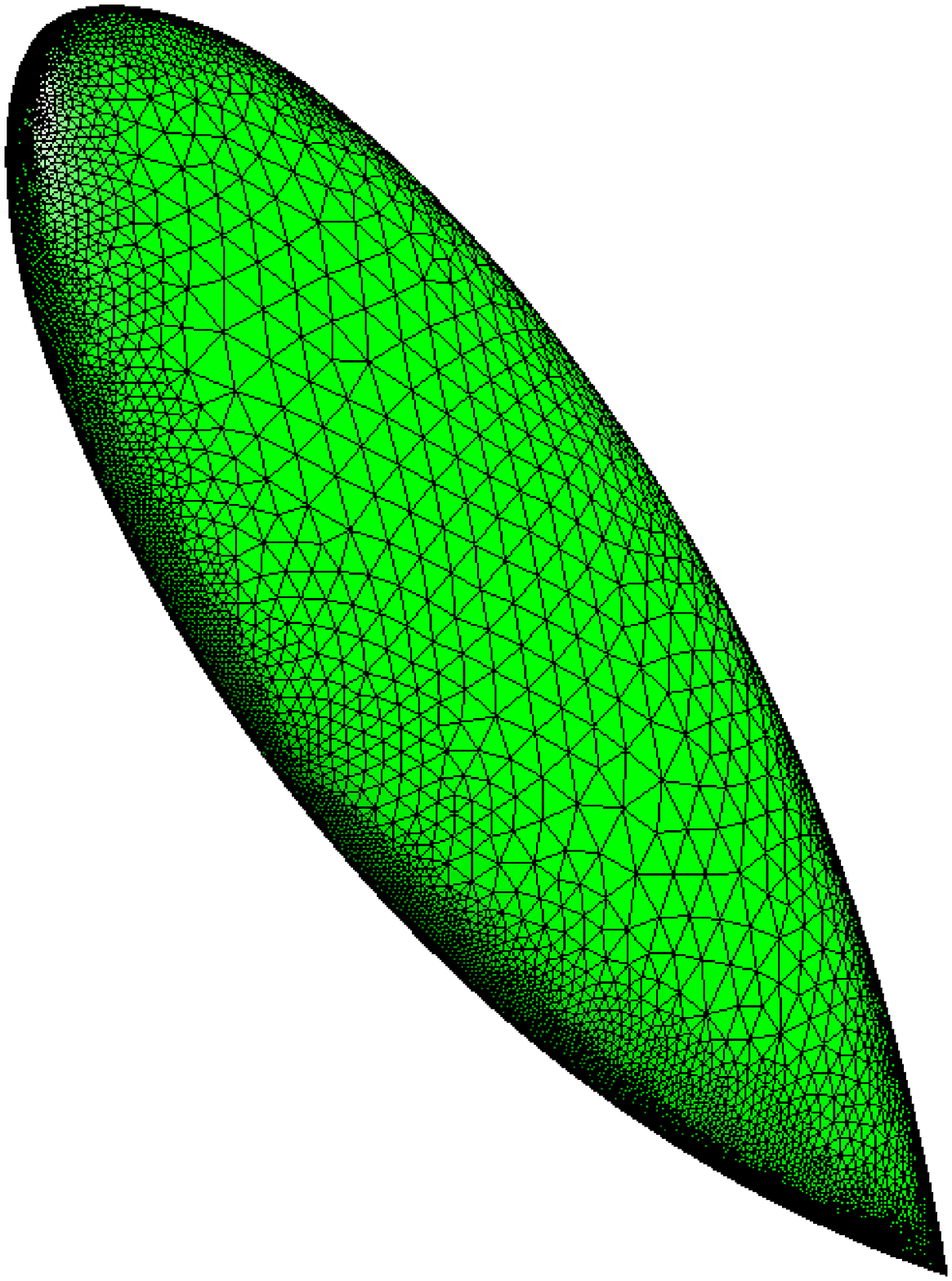}
\caption{Shapes used in the simulations. From top to bottom and from left to right: the sphere $S^{2}$ (131072 elements), a crankshaft (169600 elements), the Stanford bunny (157672 elements) and the NASA almond (116448 elements).}
\label{fig:Shapes}
\end{figure}

\subsection{Comparison of the memory usage}
In this first subsection we compare the memory required to store the $\mathcal{H}$-matrices $S^{\mathcal{H}}$, $\hat{S}^{\mathcal{H}}$, $D^{\mathcal{H}}$ and $\hat{D}^{\mathcal{H}}$. In all experiments, ACA+ with a tolerance of $\num{1e-5}$ was the method of low-rank approximation and the standard $\mathcal{H}-\text{lib}^{\text{pro}}$ weak geometric admissibilty criterion, i.e. criterion \ref{eq:admiss2} with axis-aligned bounding boxes used for distance and diameter computations, was used with separation parameter $\eta=2$. The clustering was done based on the automatic built-in partition strategy present in $\mathcal{H}-\text{lib}^{\text{pro}}$. Rather than using the absolute wavenumber, the dimensionless wavenumber is used, the product of the wavenumber $\kappa$ with the diameter of the shape. After ACA+ compression, SVD-based recompression with tolerance $\num{1e-2}$ times the ACA+ tolerance was used.

\begin{figure}[ht]
\center
\begin{tikzpicture}
\begin{axis}[width=.45\linewidth,
  ylabel near ticks,
  xlabel=$\kappa\cdot\text{diam}(\Gamma)$,
  ylabel=memory (GB),
  ylabel style={font=\small},
  legend style={font=\tiny},
  grid,
  grid style={dotted},
  legend pos=north west]
\addplot [color=blue,mark=*] table {memory_not_extracted_sphere.txt};
\addlegendentry{Sphere}
\addplot [color=blue,mark=*,dashed,forget plot] table {memory_extracted_sphere.txt};
\addplot [color=red,mark=*] table {memory_not_extracted_bunny.txt};
\addlegendentry{Bunny}
\addplot [color=red,mark=*,dashed] table {memory_extracted_bunny.txt};
\end{axis}
\end{tikzpicture}
\begin{tikzpicture}
\begin{axis}[width=.45\linewidth,
  xlabel=$\kappa\cdot\text{diam}(\Gamma)$,
  grid,
  grid style={dotted},
  legend style={font=\tiny},
  legend pos=north west]
\addplot+  table {memory_not_extracted_crankShaft.txt};
\addlegendentry{Crank}
\addplot [color=blue,mark=*,dashed,forget plot] table {memory_extracted_crankShaft.txt};
\addplot [color=red,mark=*] table {memory_not_extracted_almond.txt};
\addlegendentry{Almond}
\addplot [color=red,dashed,mark=*] table {memory_extracted_almond.txt};
\end{axis}
\end{tikzpicture}
\caption{Comparison of the memory usage for the $\mathcal{H}$-matrices $S^{\mathcal{H}}$ (full) and $\hat{S}^{\mathcal{H}}$ (dashed) for the 4 meshes given in fig. \ref{fig:Shapes} as a function of the dimensionless wave number $\kappa\cdot\text{diam}(\Gamma)$.}

\end{figure}

\begin{figure}[ht]
\begin{center}
\begin{tikzpicture}
\begin{axis}[width=.45\linewidth,
  ylabel near ticks,
  xlabel=$\kappa\cdot\text{diam}(\Gamma)$,
  ylabel=memory (GB),
  ylabel style={font=\small},
  legend style={font=\tiny},
  grid,
  grid style={dotted},
  legend pos=north west]
\addplot [color=blue,mark=*] table {memory_not_extracted_sphere_DLP.txt};
\addlegendentry{Sphere}
\addplot [color=blue,mark=*,dashed,forget plot] table {memory_extracted_sphere_DLP.txt};
\addplot [color=red,mark=*] table {memory_not_extracted_bunny_DLP.txt};
\addlegendentry{Bunny}
\addplot [color=red,mark=*,dashed] table {memory_extracted_bunny_DLP.txt};
\end{axis}
\end{tikzpicture}
\begin{tikzpicture}
\begin{axis}[width=.45\linewidth,
  xlabel=$\kappa\cdot\text{diam}(\Gamma)$,
  grid,
  grid style={dotted},
  legend style={font=\tiny},
  legend pos=north west]
\addplot+  table {memory_not_extracted_crankShaft_DLP.txt};
\addlegendentry{Crank}
\addplot [color=blue,mark=*,dashed,forget plot] table {memory_extracted_crankShaft_DLP.txt};
\addplot [color=red,mark=*] table {memory_not_extracted_almond_DLP.txt};
\addlegendentry{Almond}
\addplot [color=red,dashed,mark=*] table {memory_extracted_almond_DLP.txt};
\end{axis}
\end{tikzpicture}
\caption{Comparison of the memory usage for the $\mathcal{H}$-matrices $D^{\mathcal{H}}$ (full) and $\hat{D}^{\mathcal{H}}$ (dashed) for the 4 meshes given in fig. \ref{fig:Shapes} as a function of the dimensionless wave number $\kappa\cdot\text{diam}(\Gamma)$.}
\end{center}
\end{figure}
For both the single- and double-layer potential the extracted method is shown to be more resilient to increasing rank, staying of nearly constant memory over the wave number, while the classical method's memory usage increases linearly over the wave number. In particular the memory usage of the single-layer potential is reduced greatly by employing frequency extraction. In the case of the sphere the extracted SLP matrix $\hat{\mathcal{S}}^{\mathcal{H}}$ requires less than half the memory as the classical $\mathcal{S}^{\mathcal{H}}$. One notable exception is the NASA Almond, which is not improved upon greatly. This is due to the fact that the NASA Almond is in some respect already optimal for the BEM method: it is an elongated, anisotropic mesh, that is coarse and nearly flat in the parts of the shape that are close to each other (the two broad sides of the almond), meaning that the discretised operator's nearfield is very small compared to its farfield. Nevertheless our method is still able to improve upon the memory usage, if somewhat less than for more general shapes. 

\subsection{Comparison of computation times}
In this second subsection we compare the time required to construct the $\mathcal{H}$-matrices $S^{\mathcal{H}}$, $\hat{S}^{\mathcal{H}}$, $D^{\mathcal{H}}$ and $\hat{D}^{\mathcal{H}}$. In all experiments, ACA+ with a tolerance of $\num{1e-5}$ was the method of low-rank approximation and the standard $\mathcal{H}-\text{lib}^{\text{pro}}$ weak geometric admissibility criterion was used. The clustering was done based on the automatic built-in partition strategy present in $\mathcal{H}-\text{lib}^{\text{pro}}$. Rather than using the absolute wavenumber, the dimensionless wavenumber is used, the product of the wavenumber $\kappa$ with the diameter of the shape.

\begin{figure}[ht]
\begin{center}
\begin{tikzpicture}
\begin{axis}[width=.45\linewidth,
  ylabel near ticks,
  xlabel=$\kappa\cdot\text{diam}(\Gamma)$,
  ylabel=time (s),
  ylabel style={font=\small},
  legend style={font=\tiny},
  grid,
  grid style={dotted},
  legend pos=north west]
\addplot [color=blue,mark=*] table {timings_not_extracted_sphere.txt};
\addlegendentry{Sphere}
\addplot [color=blue,mark=*,dashed,forget plot] table {timings_extracted_sphere.txt};
\addplot [color=red,mark=*] table {timings_not_extracted_bunny.txt};
\addlegendentry{Bunny}
\addplot [color=red,mark=*,dashed] table {timings_extracted_bunny.txt};
\end{axis}
\end{tikzpicture}
\begin{tikzpicture}
\begin{axis}[width=.45\linewidth,
  xlabel=$\kappa\cdot\text{diam}(\Gamma)$,
  grid,
  grid style={dotted},
  legend style={font=\tiny},
  legend pos=north west]
\addplot+  table {timings_not_extracted_crankShaft.txt};
\addlegendentry{Crank}
\addplot [color=blue,mark=*,dashed,forget plot] table {timings_extracted_crankShaft.txt};
\addplot [color=red,mark=*] table {timings_not_extracted_almond.txt};
\addlegendentry{Almond}
\addplot [color=red,dashed,mark=*] table {timings_extracted_almond.txt};
\end{axis}
\end{tikzpicture}
\caption{Comparison of the construction time for the $\mathcal{H}$-matrices $S^{\mathcal{H}}$ (full) and $\hat{S}^{\mathcal{H}}$ (dashed) for the 4 meshes given in fig. \ref{fig:Shapes} as a function of the dimensionless wave number $\kappa\cdot\text{diam}(\Gamma)$.}
\end{center}
\end{figure}

\begin{figure}[ht]
\begin{center}
\begin{tikzpicture}
\begin{axis}[width=.45\linewidth,
  ylabel near ticks,
  xlabel=$\kappa\cdot\text{diam}(\Gamma)$,
  ylabel=time (s),
  ylabel style={font=\small},
  legend style={font=\tiny},
  grid,
  grid style={dotted},
  legend pos=north west]
\addplot [color=blue,mark=*] table {timings_not_extracted_sphere_DLP.txt};
\addlegendentry{Sphere}
\addplot [color=blue,mark=*,dashed,forget plot] table {timings_extracted_sphere_DLP.txt};
\addplot [color=red,mark=*] table {timings_not_extracted_bunny_DLP.txt};
\addlegendentry{Bunny}
\addplot [color=red,mark=*,dashed] table {timings_extracted_bunny_DLP.txt};
\end{axis}
\end{tikzpicture}
\begin{tikzpicture}
\begin{axis}[width=.45\linewidth,
  xlabel=$\kappa\cdot\text{diam}(\Gamma)$,
  grid,
  grid style={dotted},
  legend style={font=\tiny},
  legend pos=north west]
\addplot+  table {timings_not_extracted_crankShaft_DLP.txt};
\addlegendentry{Crank}
\addplot [color=blue,mark=*,dashed,forget plot] table {timings_extracted_crankShaft_DLP.txt};
\addplot [color=red,mark=*] table {timings_not_extracted_almond_DLP.txt};
\addlegendentry{Almond}
\addplot [color=red,dashed,mark=*] table {timings_extracted_almond_DLP.txt};
\end{axis}
\end{tikzpicture}
\caption{Comparison of the construction time for the $\mathcal{H}$-matrices $D^{\mathcal{H}}$ (full) and $\hat{D}^{\mathcal{H}}$ (dashed) for the 4 meshes given in fig. \ref{fig:Shapes} as a function of the dimensionless wave number $\kappa\cdot\text{diam}(\Gamma)$.}
\end{center}
\end{figure}
These timings were obtained by averaging over 5 experiments for each wave number. For both the single- and double-layer potential the timings align well with the memory usage. Even though the extracted method techically performs more FLOP per bit, our method still outperforms the classical method since less bits are needed. This means we compte far less matrix entries, ech of which corresponds to an expensive double integral of the form \ref{eq:matrixElement}. Again the construction time for single-layer potential is reduced the most by employing frequency extraction. In the case of the sphere the extracted SLP matrix $\hat{\mathcal{S}}^{\mathcal{H}}$ requires less than half the time to construct compared to the classical $\mathcal{S}^{\mathcal{H}}$. The same remark as in the previous section holds for the NASA Almond.
\subsection{Asymptotics in the grid size}
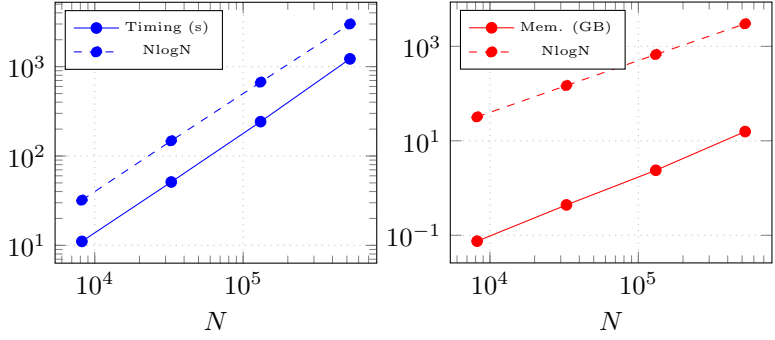
\begin{figure}[ht]\label{fig:asymptotics}
\begin{center}
\begin{tikzpicture}
\begin{axis}[width=.45\linewidth,
  xmode=log,
  ymode=log,
  ylabel near ticks,
  xlabel=$N$,
  grid,
  grid style={dotted},
  legend style={font=\tiny},
  legend pos=north west]
\addplot [color=blue,mark=*] table {sphereAsymptoticTimings.txt};
\addlegendentry{Timing (s)}
\addplot [color=blue,mark=*,dashed] table {NlogNsphere.txt};
\addlegendentry{NlogN}
\end{axis}
\end{tikzpicture}
\begin{tikzpicture}
\begin{axis}[width=.45\linewidth,
  xmode=log,
  ymode=log,
  xlabel=$N$,
  grid,
  grid style={dotted},
  legend style={font=\tiny},
  legend pos=north west]
\addplot [color=red,mark=*] table {sphereAsymptoticMem.txt};
\addlegendentry{Mem. (GB)}
\addplot [color=red,mark=*,dashed] table {NlogNsphere.txt};
\addlegendentry{NlogN}
\end{axis}
\end{tikzpicture}
\caption{Loglog plot of time and memory consumption for the construction of the extracted SLP matrix for iteratively refined sphere. Number of DOFs are $N=2^{13},2^{15},2^{17},2^{19}$. In all cases $\kappa h$ was kept constant at $.8$, with $h$ the radius of the largest triangular mesh element.}
\end{center}
\end{figure}

\begin{figure}[ht]\label{fig:ranksAnd perDofs}
\begin{center}
\begin{tikzpicture}
\begin{axis}[width=.45\linewidth,
	xmode=log,
  ylabel near ticks,
  xlabel=$N$,
  grid,
  grid style={dotted},
  legend style={font=\tiny},
  legend pos=north west]
\addplot [color=blue,mark=*] table {meanRanksSphere.txt};
\addlegendentry{Mean rank};a
\end{axis}
\end{tikzpicture}
\begin{tikzpicture}
\begin{axis}[width=.45\linewidth,
  xmode=log,
  xlabel=$N$,
  grid,
  grid style={dotted},
  legend style={font=\tiny},
  legend pos=north west]
\addplot [color=red,mark=*] table {timePerDofSphere.txt};
\addlegendentry{Mem./DOF (kB)}
\end{axis}
\end{tikzpicture}
\caption{Plot of mean rank (left) and plot of Memory per DOF of the extracted SLP matrix for iteratively refined sphere. Number of DOFs are $N=2^{13},2^{15},2^{17},2^{19}$. In all cases $\kappa h$ was kept constant at $.8$, with $h$ the diameter of the largest triangular mesh element.}
\end{center}
\end{figure}
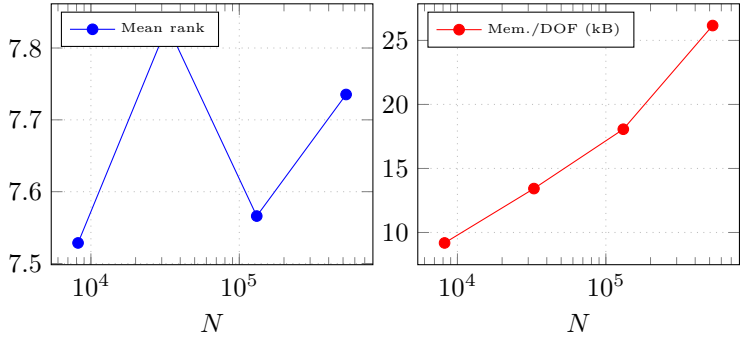

In this subsection we briefly look at the asymptotics of our extracted operator construction as a function of grid size. To this end, we discretised the sphere at 4 increasing levels of fineness. The most coarse triangulisation has $N=2^{13}$ DOFs, whereas the finest has $N=2^{19}$. The wavenumber was at every level chosen such that $\kappa h=0.8$, where $h$ is the radius of the largest element. In this way, each simulation corresponds to the construction of the total extracted SLP far-field in the high-frequency regime. The wave number dependence of the complexity and memory usage was already shown to be (weakly) linear in the previous sections. The results of the asymptotics with respect to $N$ are reported in figure \ref{fig:asymptotics}. We see clearly that the complexity is $\mathcal{O}(N\log N)$ for both memory and time. We have also investigated the mean ranks and the memory per DOF for the same resolutions of the sphere. Our findings are reported in figure \ref{fig:ranksAnd perDofs}. Note that the mean rank is nearly constant over increasing $N$, and that the memory per DOF follows an $N\log N$ type curve. We have only reported the SLP results, the DLP results are similar. Compare with the recent \cite{BormHybrid}, where the ranks still seem to increase linearly and the increase in memory per DOF is more noticable for increasing grid sizes.

\section{Numerical results: compact representation}
In this section we give the numerical results for the compact representation, for both $\hat{S}(\kappa)$ and $\hat{D}(\kappa)$, as approximated by $\hat{S}^{\mathcal{H}}(\kappa)$ and $\hat{D}^\mathcal{H}(\kappa)$. In both cases we compute a compact representation for the range $\kappa\in[a,b]$ and compare the resulting timings and memory with that of non-extracted $\mathcal{H}$-matrices at the ends of the range. The compact representation is constructed and tested only on a subset of the admissible blocks $\text{Adm}\subset\mathcal{T}_{I\times J}$, due to memory constraints affecting the performance. Improvement of memory management is a subject for future research. Let $A\subset \text{Adm}$ be this subset of the admissible blocks. Then, for $M\in\mathbb{C}^{I\times J}$, let 
\begin{align*}
\|M\|_{F,A}^2&:=\sum_{\tset\times\sset\in A}\|M_{|\tset\times\sset}\|_{F,m}^2\\
\|M\|_{\max,A}&:=\max_{\tset\times\sset\in A}\|M_{|\tset\times\sset}\|_{\max,m}
\end{align*}
with $\|.\|_{F,m}$ as in definition \ref{def:errEst} and $\|.\|_{\max,m}$ defined analogously. We compute estimates for the $\|.\|_F$ and $\|.\|_{\max}$ errors of the compact representation as
\begin{align*}
\textbf{Err}_{F}^2&:=\sum_{\kappa\in S^{\ast}}\|\hat{B}^{\otimes}(\kappa)-\hat{B}(\kappa)\|_{F,A}^2/\sum_{\kappa\in S^{\ast}}\|\hat{B}(\kappa)\|_{F,A}^2\\
\textbf{Err}_{\infty}&:=\max_{\kappa\in S^{\ast}}\|\hat{B}^{\otimes}(\kappa)-\hat{B}(\kappa)\|_{\max,A}/\max_{\kappa\in S^{\ast}}\|\hat{B}(\kappa)\|_{\max,A}
\end{align*}
for $S^{\ast}\subset [a,b]$. The set $S^{\ast}$ is chosen such that if $S$ is the $\textbf{aaa}$ sample set (the data for Algorithm~\ref{alg:linear AAA-ACA}), then $S\cap S^{\ast}=\emptyset$. 
We abbreviate $\text{Mem}(a):=\text{Mem}(B^{\mathcal{H}}(a))$ ($B=S,D$) and similarly $\text{Mem}(b):=\text{Mem}(B^{\mathcal{H}}(b))$. In the same way we write $\text{t}(a):=\text{time}(B^{\mathcal{H}}(a))$ and $\text{t}(b):=\text{time}(B^{\mathcal{H}}(b))$. The set $S$ for the $\textbf{aaa}$ routine was in all cases set to the set of $16$ Chebyshev nodes, while the maximal degree of the rational approximation was set to $8$.

\subsection{Compact representation timings and error}

We start off by approximating the extracted SLP $\hat{S}$ in the range $[a,b]=[10,100]$. For all shapes the tolerance $\text{tol}$ for the Algorithm~\ref{alg:linear AAA-ACA} is set to $\num{1e-4}$, while the tolerance for the internal $\textbf{aca}$ routine is set to $\num{1e-5}$, as in the previous experiments. This gives us table~\ref{table:SLPREPR}. 
\begin{table}\label{table:SLPREPR}
\center
 \begin{tabular}{|c|c|c|c|c|c|c|c|c|} 
 \hline
\textbf{Shape} & \textbf{Mem} & \textbf{time} & $\frac{\textbf{Mem}}{\textbf{Mem}(a)}$&$\frac{\textbf{Mem}}{\textbf{Mem}(b)}$& $\frac{t}{t(a)}$&$\frac{t}{t(b)}$&$\textbf{Err}_F$&$\textbf{Err}_{\infty}$ \\ 
 \hline
 Sphere & 13.14GB & 668.88s & 2.52 & 1.87 & 4.08 & 1.88 &1.01e-5 & 1.01e-4\\ 
 \hline
 Crank& 14.5GB&1176.72s & 2.33 & 1.54 & 3.65 & 2.25 & 1.15e-5 & 6.26e-5 \\ 
 \hline
 Bunny& 19.09GB&1122.84s & 2.40 & 1.52 & 2.42 & 1.82 &1.85e-5 & 8.51e-5 \\ 
 \hline
 Almond&11.16GB &919.86s & 2.57 & 1.96 & 3.39 & 2.81 & 1.55e-5& 7.44e-5\\ 
 \hline
\end{tabular}
\caption{Memory and time required for the construction of the compact representation of the SLP operator the four used shapes, together with two measures of the reconstruction error. Memory and time are calculated absolutely as well as relative to the required memory and time for the construction of $\mathcal{H}$-matrices at the extremes of the wave number range $[a,b]$.}
\end{table}
For the double layer potential we repeated the same experiments. Setting again the tolerance $\text{tol}$ from Algorithm~\ref{alg:linear AAA-ACA} to $\num{1e-4}$ and the tolerance for its internal $\textbf{aca}$ routine to $\num{1e-5}$, we now get figure~\ref{table:DLPREPR}.
\begin{table}\label{table:DLPREPR}
\center
 \begin{tabular}{|c|c|c|c|c|c|c|c|c|} 
 \hline
\textbf{Shape} & \textbf{Mem} & \textbf{time} & $\frac{\textbf{Mem}}{\textbf{Mem}(a)}$&$\frac{\textbf{Mem}}{\textbf{Mem}(b)}$& $\frac{t}{t(a)}$&$\frac{t}{t(b)}$&$\textbf{Err}_F$&$\textbf{Err}_{\infty}$ \\ [0.5ex] 
 \hline
 Sphere& 10.17GB & 819s & 3.09 & 1.56 & 4.61 & 2.22 &1.71e-5& 8.08e-5\\ 
 \hline
 Crank &20.22GB&2800.02s  & 2.97 & 2.29 & 5.56 & 4.28 & 1.18e-5& 1.e-4\\ 
 \hline
 Bunny &19.43GB & 2637s  & 3.07 & 2.27 & 5.68 & 4.28 & 2.14e-5&1.01e-4 \\ 
 \hline
 Almond&11.60GB & 1826s & 3.24 & 2.71 & 6.73 & 5.57 &1.09e-5 & 9.56e-5\\ 
 \hline
\end{tabular}
\caption{Memory and time required for the construction of the compact representation of the DLP operator the four used shapes, together with two measures of the reconstruction error. Memory and time are calculated absolutely as well as relative to the required memory and time for the construction of $\mathcal{H}$-matrices at the extremes of the wave number range $[a,b]$.}
\end{table}
These experiments show that the construction of our compact representation can be done in a reasonable time and memory complexity, when compared to the original construction of the non-extracted $\mathcal{H}$-matrices. For the single layer potential the construction of the compact representation can always be done in less than 4 times the time it takes to construct a single corresponding $\Hmat$, while for the double layer potential (with the exception of the NASA Almond), this can be done in less than 6 times the construction time for the corresponding $\Hmat$. We have shown that the error is well within bounds and the memory usage is proportionally low as well.

\subsection{Reconstruction timings}
In this section we compare the computation time of the reconstruction of the far-field blocks from the compact representation ($\textbf{time}_R$) to the computation time of the construction of those same blocks in the classical and extracted $\mathcal{H}$-matrix approach ($\textbf{time}_O$ and $\textbf{time}_E$), at the end of the range $[a,b]$. We have chosen $\kappa=90$ as this is a high-frequency wave number that does not lie in the interpolation set that we have used (the $16$ Chebyshev nodes in $[10,100]$). The results are given in tables \ref{table:SLPREC} and \ref{table:DLPREC}. The near-field is not included in this computation.
\begin{table}\label{table:SLPREC}
\center
 \begin{tabular}{|c|c|c|c|} 
 \hline
\textbf{Shape} & $\textbf{time}_R$ & $\textbf{time}_E$ & $\textbf{time}_O$ \\ [0.5ex]  
 \hline
 Sphere & 10.50 & 147.36 & 374.28 \\ 
 \hline
 Crank & 22.92 & 249.06 & 456.24\\ 
 \hline
 Bunny & 17.88 & 262.32 & 430.86\\ 
 \hline
 Almond & 12.48 & 180.36 & 229.68\\ 
 \hline
\end{tabular}
\caption{Reconstruction timings for the SLP.}
\end{table}

\begin{table}\label{table:DLPREC}
\center
 \begin{tabular}{|c|c|c|c|c|c|} 
 \hline
\textbf{Shape} & $\textbf{time}_R$ & $\textbf{time}_E$ & $\textbf{time}_O$ & $\frac{\textbf{time}_E}{\textbf{time}_R}$& $\frac{\textbf{time}_O}{\textbf{time}_R}$ \\ [0.5ex] 
 \hline
 Sphere & 12.36 & 164.19 & 371.64 &13.28&30.07\\ 
 \hline
 Crank & 69.18 & 439.32 & 561.02 & 6.35&8.11\\ 
 \hline
 Bunny & 42.48 & 480.84 & 518.94 & 11.32&12.22 \\ 
 \hline
 Almond& 24.84 & 257.34 & 298.56 & 10.36&12.02\\ 
 \hline
\end{tabular}
\caption{Reconstruction timings for the DLP.}
\end{table}
As can be seen, with the exception of the DLP for the NETGEN Crankshaft (our most complex case due to its many interacting flat clusters, see \cite{HCA}) our method provides an acceleration with a factor $10$, a remarkable gain. Even for the NETGEN Crankshaft, we have a gain of a factor $6$. This, in combination with the previous subsection shows heuristically that as soon as more than $4$ wave numbers in the case of the SLP or $6$ in case of the DLP are to be considered, our compact representation method will provide significant gains over our extracted method, which itself often provides significant gains over the classical method. The main reason for the large speedup is that through the use of the compact representation we avoid having to calculate the double integrals of the form in equation \ref{eq:matrixElement}.
\section{Conclusion and future work}
We have shown that frequency extraction is a powerful tool for the construction of high-frequency $\mathcal{H}$-matrices and for the construction of a compact representation of their wave number dependence. We have demonstrated theoretically and experimentally that our frequency extraction results in a nearly constant memory and time complexity over a large range of wave numbers. We have also demonstrated that a compact wave number dependence can be constructed cheaply and exploited for the fast reconstruction of extracted $\mathcal{H}$-matrices. We have shown that the gain of this reconstruction over the extracted method (and hence the classical method) is enormous.

In future work we will provide a fast matrix-vector product of an extracted BEM-matrix with a vector, based on the Hadamard product structure. We will also implement a more sophisticated algorithm with customised memory management. In addition, we will use CORK and set-valued AAA (see \cite{CORK} and \cite{KarlSVAAA}) to further reduce the complexity of the compact representation and implement a frequency sweep.
\section*{Acknowledgements}
This work was supported in part by KU Leuven IF project C14/15/055 and FWO Research Foundation Flanders G0B7818N.
We would like to thank the referees for their thorough comments and suggestions.
\appendix
\section{Separability of the extracted kernel}\label{appendix:A}
In this appendix we will show Theorem~\ref{thm:LRindep}. Given some vector $\mathbf{r}\in\mathbb{R}^3$, the symbols $r$ and $\hat{\mathbf{r}}$ will refer to its norm and unit direction respectively. The radial functions $j_n(r)$ and $y_n(r)$ are the \emph{spherical Bessel functions of the first and second kind} respectively, while $h_n^{(1)}(r)=j_n(r)+\imath y_n(r)$ are the spherical Hankel functions of the first kind, see \cite{Abramowitz}. We fix the following conventions:
\begin{definition}[Legendre Polynomials and Spherical Harmonics]
The \emph{Legendre polynomials} $\{P_n\}_{n\in\mathbb{N}}$ and the \emph{Associated Legendre Polynomials} $\{P_n^m\}_{n,m\in\mathbb{N}}$ on $[-1,1]$ are defined for $n\geq m\geq 0$ by
$$P_n(x)=\frac{1}{2^nn!}\frac{\D^n}{\D x^n}(x^2-1)^n,\qquad P_n^m(x)=(-1)^{m}(1-x)^{m/2}\frac{\D^m}{\D x^m}P_n(x)$$
and $P_n^{-m}=(-1)^m\frac{(n-m)!}{(n+m)!}\cdot P_{n}^m$.\newline The \emph{Spherical Harmonics} $\{Y_n^m|n\in\mathbb{N},m=-n,\ldots,n\}$ on $[0,\pi]\times[0,2\pi)$ are defined by
$$Y_n^m(\theta,\phi)=\sqrt{\frac{2n+1}{4\pi}\frac{(n-m)!}{(n+m)!}}P_n^{m}(\cos{\theta})e^{\imath m\phi}.$$
When $\mathbf{r}$ satisfies $\mathbf{r}=(r,\theta,\phi)$ in spherical coordinates we will use the notation $Y_n^m(\hat{\mathbf{r}}):=Y_n^m(\theta,\phi)$. In addition we will abbreviate $\int_0^{\pi}\sin{\theta}\D\theta\int_{0}^{2\pi}\D\phi=:\int_{S^2}\D\Omega(\mathbf{s})$.
\end{definition}
The above conventions imply that $\|P_n(x)\|_{\infty}=\sup_{x\in[-1,1]}|P(x)|=1$ and $\int_{S^2}Y_{n}^{m}Y_{n'}^{m'}\D\Omega(\mathbf{s})=\delta_{nn'}\delta_{mm'}$. We also have the well-known \emph{Addition Theorem}:
\begin{equation}\label{eq:AdditionThmLegendre}
P_n(\hat{\mathbf{r}}_1\cdot\hat{\mathbf{r}}_2)=\frac{4\pi}{2n+1}\sum_{m=-n}^nY^{-m}_n(\hat{\mathbf{r}}_1)Y^{m}_n(\hat{\mathbf{r}}_2)
\end{equation}
for all $n\in\mathbb{N}$ and $\mathbf{r}_1,\mathbf{r}_2\in\mathbb{R}^3$ (see e.g. \cite{Ferrers}). Also, from e.g. \cite[\S4,2.19(b)]{EuclidFourier} we have
\begin{equation}\label{eq:normBoundY}
\|Y_n^m\|_{\infty}=\sup_{\mathbf{s}\in S^2}|Y_{n}^m(\hat{\mathbf{s}})|\leq \sqrt{\frac{n+1}{2\pi}}.
\end{equation}
Our first observation is that local separable expansions, at the level of individual DOF-DOF interactions can be glued together:
\begin{lemma}[Gluing lemma]\label{lem:gluedSep}
Suppose $\tset=\cup_i\tau_i$, $\sset=\cup_j\sigma_j$ are two clusters in $\Gamma$, with $\forall i_1\neq i_2: \tau_{i_1}\cap\tau_{i_2}=\emptyset$ and $\forall j_1\neq j_2: \sigma_{j_1}\cap\sigma_{j_2}=\emptyset$. Suppose the kernel
$b:\tset\times\sset\to\mathbb{C}:(\mathbf{x},\mathbf{y})\mapsto b(\mathbf{x},\mathbf{y})$ satisfies
$b_{|\tau_i\times \sigma_j}(\mathbf{x},\mathbf{y})=\sum_{\nu=1}^{R} u^{(i)}_{\nu}(\mathbf{x}) v^{(j)}_{\nu}(\mathbf{y})$,
with $\forall i,j: u^{(i)}_{\nu}:\tset\to \mathbb{C}, v^{(j)}_{\nu}:\sset\to \mathbb{C}$ arbitrary functions.
Then there exist functions $u_{\nu}:\tset\to \mathbb{C}, v_{\nu}:\sset\to \mathbb{C}$, $\nu=1,\ldots,R$ such that $b$ can be written as 
$$b(\mathbf{x},\mathbf{y})=\sum_{\nu=1}^{R} u_{\nu}(\mathbf{x})v_{\nu}(\mathbf{y})$$
\end{lemma}
\begin{proof}
The proof of this is simple. Let $\chi_{i}$ and $\chi_j$ denote the indicator functions of $\tau_i$ and $\sigma_j$ respectively. Then, clearly,
\begin{equation*}
b(\mathbf{x},\mathbf{y})=\sum_{ij}\sum_{\nu=1}^{R}\chi_{i}(\mathbf{x})u^{(i)}_{\nu}(\mathbf{x})\chi_{j}(\mathbf{y}) v^{(j)}_{\nu}(\mathbf{x})=\sum_{\nu=1}^{R}\underbrace{\left(\sum_{i}(\chi_{i} u^{(i)}_{\nu})(\mathbf{x})\right)}_{:= u_{\nu}}\cdot\underbrace{\left(\sum_{j}(\chi_{j} v^{(j)}_{\nu})(\mathbf{y})\right)}_{:= v_{\nu}}
\end{equation*}
which concludes the proof.
\end{proof}
Next we will consider three different forms of addition theorems, each of which we use in the final proof. Theorems \ref{thm:TransGreen} and \ref{thm:regTrans} are specialised restatements of the same well-known addition theorem for solutions to the Helmholtz eqution, which can be found in \cite[D.14]{chewWaves}. In particular we derive specialised truncation error bounds. Theorem \ref{thm:TransSolid} is also a well-known result, though our proof of its convergence rate seems to be new.
\begin{theorem}[Solid Harmonic Translation Theory]\label{thm:TransSolid}
Let $n\in\mathbb{N}$, and $m\in\{-n,\ldots,n\}$. Let $I_n^{m}(\mathbf{r}):=Y_n^m(\hat{\mathbf{r}})/r^{n+1}$ and $R_n^{m}(\mathbf{r}):=r^nY_n^m(\hat{\mathbf{r}})$ denote the (scaled) irregular and regular solid harmonics. Let $\D\mathbf{r}\in\mathbb{R}^3$ such that $\rho:=\D r/r<1$. Then
$$I_n^{m}(\mathbf{r}+\D\mathbf{r})=2\sqrt{\pi}\sum_{\nu=0}^{\infty}\sum_{\mu=-\nu}^{\nu}\frac{(-1)^{\nu+\mu}c_{\nu\mu}}{\sqrt{2\nu+1}}R_{\nu}^{\mu}(\D\mathbf{r})I_{n+\nu}^{m-\mu}(\mathbf{r})$$
with 
$$c_{\nu\mu}:=\left(\frac{(2n+1)}{(2n+2\nu+1)}\binom{n+m+\nu-\mu}{\nu-\mu}\binom{n-m+\nu+\mu}{\nu+\mu}\right)^{1/2}.$$
For any $0<\delta\leq 1$ there is a degree $2n+2$ monic polynomial $q$ such that the truncation error of the above expansion is given by
\begin{equation}\label{eq:truncationErr}
\epsilon_p:=\left|2\sqrt{\pi}\hspace*{-5pt}\sum_{\nu=p+1}^{\infty}\sum_{\mu=-\nu}^{\nu}\frac{(-1)^{\nu+\mu}c_{\nu\mu}}{\sqrt{2\nu+1}}R_{\nu}^{\mu}(\D\mathbf{r})I_{n+\nu}^{m-\mu}(\mathbf{r})\right|<\sqrt{\frac{3q\left(2N_{\delta}^{(n)}\right)}{4\pi(2n)!}}\frac{(1+\delta)^{p+1-N_{\delta}^{(n)}}\hspace*{-5pt}\rho^{p+1}}{(1-\rho)r^{n+1}}
\end{equation}
for all $p\geq N_{\delta}^{(n)}:=(n+1)/\delta$.
\end{theorem}
\begin{proof}
This expansion can be found in \cite{Piecuch}. All that is left to show is its convergence rate. Suppose $n>0$ and $m \geq 0$. Then, naturally
$$\frac{(2n+1)}{(2n+2\nu+1)}\leq 3\frac{n+m}{n+m+\nu-\mu}$$
and so, by the AM-QM inequality and the Rothe-Hagen identity, we have
\begin{align*}
\sum_{\mathclap{\mu=-\nu}}^{\nu}c_{\nu\mu}&\leq\left(3(2\nu+1)\sum_{\mathclap{\mu=-\nu}}^{\nu}\left(\frac{n+m}{n+m+\nu-\mu}\right)\binom{n+m+\nu-\mu}{\nu-\mu}\binom{n-m+\nu+\mu}{\nu+\mu}\right)^{1/2}\\
&\leq\left(3(2\nu+1)\binom{2\nu+2n}{2\nu}\right)^{1/2}
\end{align*}
The case $m\leq 0$ is proved similarly. The case $n=0$ is trivial. By equation \ref{eq:normBoundY} we get
\begin{align*}
\epsilon_p&\leq \frac{\sqrt{3}}{r^{n+1}\sqrt{\pi}}\sum_{\nu=p+1}^{\infty}\sqrt{(\nu+1)(\nu+n+1)}\binom{2\nu+2n}{2\nu}^{1/2}\rho^{\nu}=:\frac{\sqrt{3}}{\sqrt{\pi} r^{n+1}}\sum_{\nu=0}^{\infty}a_{\nu}\rho^{\nu}
\end{align*}
Note that $4a_{\nu}^2(2n)!$ is a degree $2n+2$ monic polynomial in $2\nu$ and $a_{\nu+1}/a_{\nu}\leq \frac{\nu+n+1}{\nu}$, meaning $N_{n}^{(\delta)}\leq \nu$ implies $a_{\nu+1}/a_{\nu}\leq 1+\delta$. Then  
$$\epsilon_{p}<a_{N_{n}^{(\delta)}}\rho^{N_{n}^{(\delta)}}\sum_{k=0}^{\infty}((1+\delta)\rho)^{k}.$$
The proof is then completed by truncating the left-over series in the right-hand side.
\end{proof}
Observe that that the factor $1/r^{n+1}$ does not contribute to the \emph{relative} truncation error and that, since $\rho<1/2$ holds in our application of this theorem, the restriction on convergence by $\delta$ is very weak. The following theorem we adapt from \cite{Darve}:
\begin{theorem}[Green's kernel Translation Theory]\label{thm:TransGreen}
Let $\mathbf{r},\D\mathbf{r}\in\mathbb{R}^3$ be such that $\D r/r<2\sqrt{5}/5$, $\kappa r>1$ and $\kappa\D r<1/2$. Let $\mathbf{r}':=\mathbf{r}+\D\mathbf{r}$. We have that
$$\frac{e^{\imath\kappa r'}}{\imath\kappa r'}=4\pi\sum_{n=0}^{\infty}\sum_{m=-n}^{n}h_n^{(1)}(\kappa r)Y_{n}^{m}(\hat{\mathbf{r}})j_n(\kappa \D r)Y_{n}^{-m}(\D\hat{\mathbf{r}}).$$
For all $0<\epsilon$ there is a $\tilde{C}>1$ such that if $p>\log(\tilde{C}/\epsilon)\left(\frac{1}{\log{2}}+\frac{1}{\log{\frac{\sqrt{5}}{2}}}\right)$, the truncation error $\epsilon_p$ of this expansion (defined similarly to equation \ref{eq:truncationErr}) satisfies $\epsilon_p\leq \epsilon$.
\end{theorem}
\begin{proof}
Recall that by equation \ref{eq:AdditionThmLegendre} we can write the individual terms in the above expansion as $(2n+1)j_n(\kappa\D r)h_n^{(1)}(\kappa r)P_n(\hat{\mathbf{r}}\cdot\D\hat{\mathbf{r}})$. We follow the arguments and notation from \cite{Darve}, but adapted to our case ($\kappa\D r<1$). From \cite[\S5]{Darve} we have that, when $n+\frac{1}{2}>\kappa \D r$,
$$\frac{j_{n}(\kappa\D r)}{j_{n-1}(\kappa\D r)}\leq \frac{\kappa\D r}{2n+1}G(2,n,\kappa\D r)\leq\frac{\kappa\D r}{2n+1}.$$
Similarly, whenever $n-\frac{7}{2}>\kappa r$
\begin{equation}\label{eq:besselRatio}
\frac{y_{n}(\kappa r)}{y_{n-1}(\kappa r)}\leq \frac{2n-1}{\kappa r}F(1,n,\kappa r).
\end{equation}
Then, since $\frac{(\kappa r)^2}{(2n+1)(2n-1)}G(2,n,\kappa r)F(1,n,\kappa r)^{-1}<1$ for $n-\frac{7}{2}\geq\kappa r$, we have that equation \ref{eq:besselRatio} also holds for $\frac{|h_{n}(\kappa r)|}{|h_{n-1}(\kappa r)|}$. Using $|j_{n}(\kappa\D r)|\leq 1$ in place of the bound provided in \cite[prop. 6.1]{Darve} together with $q:=\ceil{\kappa r+\frac{7}{2}}$, we get that
$$(2n+1)|j_{n}(\kappa \D r)||h_{n}(\kappa r)||P_{n}(\hat{\mathbf{r}}\cdot \D\hat{\mathbf{r}})|\leq
\begin{cases}
\frac{C(\kappa\D r)^{n-5/6}}{(2n-1)!!}& 0<n< q\\
\frac{C(\kappa\D r)^{n-5/6}}{(2n-1)!!}\left(\frac{\D r}{r}\right)^{n-q} &  q \leq n\\
\end{cases}
$$
and so, for $0<p<q$, $\epsilon_p$ is bounded by (see \cite[\S 7]{Abramowitz} for the error function $\text{erf}$) 
\begin{equation*}
\epsilon_p\leq\frac{C(\kappa\D r)^{p+1/6}}{(2p+1)!!}\sum_{n=0}^\infty \frac{1}{(2n+1)!!}=\frac{C(\kappa\D r)^{p+1/6}}{(2p+1)!!}\sqrt{\frac{\pi e}{2}}\text{erf}\left(\frac{1}{\sqrt{2}}\right).
\end{equation*}
The constant $C$ is found by requiring $|h_q(\kappa r)|<C(\kappa r)^{-5/6}$ for all $q-9/2\leq\kappa r\leq q-7/2$ and it decreases as $\kappa$ increases. For $\kappa\rightarrow 0$ this $C$ grows unbounded, but since this case is not considered here, there exists some uniform bound $\epsilon_p\leq\tilde{C}(\kappa\D r)^{p}$ for all $p<q$. Suppose that indeed $p<q$. Then for any $0<\epsilon<1$ given, $\epsilon_p<\epsilon$ when
$p>\frac{\log(\tilde{C}/\epsilon)}{\log(1/\kappa\D r)}$.
On the other hand, if convergence to precision $\epsilon$ is delayed past $p \geq q$ it must hold that $\frac{\log(\tilde{C}/\epsilon)}{\log(1/\kappa\D r)}\geq q$. For any $p\geq q$, $\epsilon_p<\epsilon$ thus holds whenever
\begin{equation}\label{eq:pBound}
p>\frac{\log(\tilde{C}/\epsilon)}{\log(1/\kappa\D r)}+\frac{\log(\tilde{C}/\epsilon)}{\log(r/\D r)} \geq q+\frac{\log(\tilde{C}/\epsilon)}{\log(r/\D r)}
\end{equation}
since we then have geometric convergence with a factor $\D r/r$. This, together with $\kappa\D r < \frac{1}{2}$ and $\frac{\D r}{r}\leq \frac{2}{\sqrt{5}}$ completes the proof.
\end{proof}
\begin{theorem}[Regular Helmholtz Translation Theory]\label{thm:regTrans} Suppose that $n\in\mathbb{N}$ and $\mathbf{r},\D\mathbf{r}\in\mathbb{R}^3$ with $\D r<r$ and $\kappa\D r<1$. Let $\mathbf{r}':=\mathbf{r}+\D\mathbf{r}$. Then, for all $m=-n,\ldots,n$:
\begin{equation}\label{eq:expansion}
j_n(\kappa r')Y_n^m(\hat{\mathbf{r}}')=\sum_{\nu=0}^{\infty}\sum_{\mu=-\nu}^{\nu}j_{\nu}(\kappa\D r)Y_{\nu}^{\mu}(\D\hat{\mathbf{r}})T_{nm}^{\nu\mu}(\kappa\mathbf{r})
\end{equation}
with
$$T_{nm}^{\nu\mu}(\kappa\mathbf{r})=\sum_{\nu'=|n\text{-}\nu|}^{n+\nu}4\pi\imath^{\nu'+\nu-n}j_{\nu'}(\kappa r)Y_{\nu'}^{m\text{-}\mu}(\hat{\mathbf{r}})\int_{S^2}Y_{n}^{m}(\mathbf{s})Y_{\nu}^{\text{-}\mu}(\mathbf{s})Y_{\nu'}^{\mu\text{-}m}(\mathbf{s})\,\D\Omega(\mathbf{s}).$$
The truncation error $\epsilon_p$ for $p \geq 1$ (defined similarly to equation \ref{eq:truncationErr}) of this expansion is given by
$$
\epsilon_{p}\leq \frac{(2n+1)\sqrt{n+1}}{(2p-1)!!}(\kappa\D r)^{p+1}.
$$
\end{theorem}
\begin{proof}
In \cite{Gumerov} this expansion is shown to hold, however, the truncation error is not shown (at least not exactly, see \cite[\S9]{Gumerov}). We aim for a concise bound, tailored to our case. As shown in \cite[\S3]{Gumerov}, due to the relation of the right-hand side integral to the \emph{Clebsch-Gordan coefficients} it holds that
$$T_{nm}^{\nu\mu}(\kappa\mathbf{r})=\sum_{\nu'=|n-\nu|}^{n+\nu}\sum_{\mu'=-\nu'}^{\nu'}4\pi\imath^{\nu'+\nu-n}j_{\nu'}(\kappa r)Y_{\nu'}^{\mu'}(\hat{\mathbf{r}})\int_{S^2}Y_{n}^{m}(\mathbf{s})Y_{\nu}^{-\mu}(\mathbf{s})Y_{\nu'}^{-\mu'}(\mathbf{s})\,\D\Omega(\mathbf{s}).$$
Using equation \ref{eq:AdditionThmLegendre} on $\sum_{\mu=\text{-}\nu}^{\nu}Y_{\nu}^{\mu}(\D\hat{\mathbf{r}})Y_{\nu}^{\text{-}\mu}(\mathbf{s})$ and $\sum_{\mu'=\text{-}\nu'}^{\nu'}Y_{\nu'}^{\mu'}(\hat{\mathbf{r}})Y_{\nu'}^{\text{-}\mu'}(\mathbf{s})$ together with the bounds $\forall\nu:\|P_{\nu}\|_{\infty},\|j_{\nu}\|_{\infty}\leq 1$ and equation \ref{eq:normBoundY} we then obtain
\begin{align*}
\epsilon_p\leq(2n+1)\sqrt{\frac{n+1}{2\pi}}\sum_{\nu=p+1}^{\infty}(2\nu+1)^2|j_{\nu}(\kappa\D r)|.
\end{align*}
From \cite{Abramowitz}, $|j_{\nu}(\kappa\D r)|\leq (\kappa\D r)^{\nu}/(2\nu+1)!!$ and $\sum_{\nu=0}^\infty 1/(2\nu+1)!!=\sqrt{\pi e/2}\,\text{erf}(1/\sqrt{2})$, so
\begin{align*}
\epsilon_p&\leq\frac{5(2n+1)}{3}\sqrt{\frac{n+1}{2\pi}}\sum_{\nu=p+1}^{\infty}\frac{(\kappa\D r)^{\nu}}{(2\nu-3)!!}\leq\frac{5\sqrt{e}(2n+1)\sqrt{n+1}(\kappa\D r)^{p+1}}{6(2p-1)!!}\,\text{erf}(1/\sqrt{2})
\end{align*}
recalling $p\geq 1$. Since $\text{erf}(1/\sqrt{2})(5\sqrt{e}/6)<1$ this concludes the proof.
\end{proof}
We are now ready to prove theorem \ref{thm:LRindep}:
\begin{proof}\label{proof:perturbation}
Note that it suffices to show that the kernel
$$\tset\times\sset\mapsto\mathbb{C}:(\mathbf{x},\mathbf{y})\mapsto\frac{\exp(\imath\kappa(\|\mathbf{x}-\mathbf{y}\|-\|\xi(\mathbf{x})-\eta(\mathbf{y})\|))}{\|\mathbf{x}-\mathbf{y}\|}$$
is separable, where $\xi(\mathbf{x})$ is defined as the center of the DOF to which $\mathbf{x}$ belongs and $\eta(\mathbf{y})$ is defined similarly. Without loss of generality we assume $\tset$ is centered around zero. Then for any DOF-DOF interaction in $\tset\times\sset$ characterised by the centers $\xi_i$, $\eta_j$ it holds that $\rho:=\|\xi_i\|/\|\eta_j\|<1/2$. Set $\mathbf{r}:=(\xi_i-\eta_j)$, $\D\mathbf{r}:=(\mathbf{x}-\xi_i)+(\eta_j-\mathbf{y})$, $\mathbf{r}':=\mathbf{r}+\D\mathbf{r}$, $\D\mathbf{x}:=(\mathbf{x}-\xi_i)$ and $\D\mathbf{y}:=(\eta_j-\mathbf{y})$. Observe that 
$$\frac{p_n(r;\kappa)}{\imath\kappa}:=r^{n+1}e^{\text{-}\imath r}h_n^{(1)}(r)=\frac{1}{\kappa\imath^{(n+1)}}\sum_{k=0}^{n}\left(\frac{i}{2\kappa}\right)^{k}\frac{(n+k)!}{k!(n-k)!}r^{n\text{-}k}$$
defines a radial polynomial $p_n(r;\kappa)$ whose non-leading coefficients decay as $\kappa\rightarrow\infty$. In addition, since $\kappa r>1$ is assumed, there is a uniform bound over $n$ and $\kappa$ for the size of these coefficients, and the norm of the leading coefficient is always 1. As a smooth function of the (asymptotically smooth) distance kernel, $p_n(\text{dist}(\mathbf{x},\mathbf{y});\kappa)$ is now separable over $\tset\times\sset$, see \cite[\S E.3]{H2Hackbusch}, with a uniform bound over $n$ and $\kappa$ on the numerical separability rank. By theorem's \ref{thm:TransSolid}, \ref{thm:TransGreen} and \ref{thm:regTrans} we have
\begin{align*}
\frac{e^{\imath\kappa(r'-r)}}{r'} &= \sum_{n=0}^{\infty}\sum_{m=\text{-}n}^m \imath \kappa e^{\text{-}\imath r}h_n^{(1)}(r)Y_n^{\text{-}m}(\mathbf{r})j_n(\kappa\D\mathbf{r})Y_{n}^{m}(\D\hat{\mathbf{r}})\\
&=\sum_{\mathclap{n,\nu_1,\nu_2=0}}^{\infty}p_n(r;\kappa)\hspace*{-5pt}\sum_{m=\text{-}n}^{n}\sum_{\mu_1=\text{-}\nu_1}^{\nu_1}\sum_{\mu_2=\text{-}\nu_2}^{\nu_2}\hspace*{-5pt}I_{n+\nu_1}^{m\text{-}\mu_1}(\xi_i)j_{\nu_2}(\kappa \D x)Y_{\nu_2}^{\text{-}\mu_2}(\D\hat{\mathbf{x}})R_{\nu_1}^{\mu_1}(\eta_j)T_{nm}^{\nu_2\mu_2}(\kappa\D\mathbf{y}).
\end{align*}
For any $\epsilon>0$ we now know by theorems \ref{thm:TransSolid}, \ref{thm:TransGreen}  and \ref{thm:regTrans} that a (relative) error of $\epsilon$ can be obtained by truncating the above sum to $n\leq p=\mathcal{O}(\log(1/\epsilon))$, $\nu_1\leq p_1=\mathcal{O}(\log{(1/\epsilon)})$ and $\nu_2\leq p_2 =\log{(\epsilon)}/\log{(\kappa\D r)}$. The total number of terms is then $\mathcal{O}(\log^6(\epsilon)/\log^2(\kappa\D r))$. Using lemma \ref{lem:gluedSep} to glue these separable representations together over $\tset\times\sset$ then concludes the proof. 
\end{proof}

\bibliographystyle{abbrv}
\bibliography{BEMApprox}
\end{document}